\documentclass[11pt,leqno,a4paper,english]{amsart}
\usepackage{amssymb,amsmath,amsfonts,amsthm,amscd}
\makeatletter
\@namedef{subjclassname@2020}{\textup{2020} Mathematics Subject Classification}
\makeatother

\usepackage{leftidx}
\usepackage{tikz}
\usepackage[utf8]{inputenc}

\usepackage{mathrsfs}

\usepackage{hyperref}
\hypersetup{
    colorlinks=true,       
    linkcolor=blue,          
    citecolor=blue,        
    filecolor=magenta,      
    urlcolor=cyan           
}

\usepackage{subfigure,graphicx}
\usepackage{xcolor}
\usepackage{xspace}
\usepackage{enumerate,enumitem}
\usepackage[all]{hypcap}

\renewcommand{\div}{ {\text{div}}}

\newcommand{\mb}[1]{\ensuremath{\mathbb{#1}}}
\newcommand{\N}{{\mb{N}}}

\newcommand{\R}{{\mb{R}}}
\newcommand{\C}{{\mb{C}}}

\newcommand{\E}{\ensuremath{\mathcal E}}
\newcommand{\G}{\ensuremath{\mathcal G}}
\renewcommand{\H}{\ensuremath{\mathcal H}}

\renewcommand{\L}{\ensuremath{\mathcal L}}

\renewcommand{\d}{\ensuremath{\partial}}

\renewcommand{\leq}{\leqslant}
\renewcommand{\geq}{\geqslant}
\newcommand{\dx}{{\, {\rm d}x}}
\newcommand{\dt}{{\, {\rm d}t}}
\usepackage{subfigure,graphicx}
\usepackage{color}
\usepackage{xspace}

\DeclareMathOperator{\supp}{Supp}

\newtheoremstyle{note}{} {}{\itshape}{-6pt}{\bf}{. --}{ }{}

\newtheorem{theorem}{Theorem}[section]

\newtheorem{lemma}[theorem]{Lemma}
\newtheorem{corollary}[theorem]{Corollary}

\newtheorem*{theo*}{Theorem}

\newcounter{theorembiss}

\newtheorem{defi}[theorem]{Definition}
\newtheorem{rema}[theorem]{Remark}

\newcommand{\Obs}{\mathrm{Obs}}

\numberwithin{equation}{section}

\subjclass[2020]{35A18, 35L05, 35Q93, 93B07, 93C20}

\title[Regional and partial observability and control of waves] 
{Regional and partial observability and control of waves} 

 \author{Belhassen Dehman}
 \address{Belhassen Dehman. D\'epartement de Math\'ematiques, Facult\'e
  des sciences de Tunis $\&$ Enit-Lamsin, Universit\'e de Tunis El Manar, 2092 El
  Manar, Tunisia. }
\email{belhassen.dehman@fst.utm.tn}

\author{Sylvain Ervedoza}
 \address{Institut de Math\'ematiques de Bordeaux, UMR 5251, Universit\'e de Bordeaux, CNRS, Bordeaux INP, F-33400 Talence, France. }
\email{sylvain.ervedoza@math.u-bordeaux.fr}

 \author{Enrique Zuazua}
 \address{Enrique Zuazua. [1] Chair for Dynamics, Control and Numerics - Alexander von Humboldt-Professorship, Department of Data Science, Friedrich-Alexander-Universit\"at Erlangen-N\"urnberg,
91058 Erlangen, Germany ,
\newline \indent\hskip 7pt
[2] Chair of Computational Mathematics, Universidad de  Deusto,
48007 Bilbao, Basque Country, Spain,
\newline \indent\hskip 7pt
[3] Departamento de Matem\'{a}ticas,
Universidad Aut\'{o}noma de Madrid,
28049 Madrid, Spain.}
\email{enrique.zuazua@fau.de}
\date{\today}
\begin{document}

\begin{abstract} 
We establish sharp regional observability results for solutions of the wave equation in a bounded domain $\Omega \subset \mathbb{R}^n$, in arbitrary spatial dimension. Assuming the waves are observed on a non-empty open subset $\omega \subset \Omega$ and that the initial data are supported in another open subset \( \mathscr{O} \subset \Omega \), we derive estimates for the energy of initial data localized in \( \mathscr{O} \), in terms of the energy measured on the observation set $(0,T) \times \omega$. This holds under a suitable geometric condition relating the time horizon $T$ and the pair of subdomains \( (\omega, \mathscr{O}) \).

Roughly speaking, this geometric condition requires that all rays of geometric optics in $\overline{\Omega}$, emanating from \( \overline{\mathscr{O}} \), must reach the observation region $(0,T) \times \omega$. Our result 
generalizes classical observability results, which recover the total energy of all solutions when the observation set $\omega$ satisfies the so-called  Geometric Control Condition (GCC)—a particular case corresponding to \( \mathscr{O} = \Omega \).

A notable feature of our approach is that it remains effective in settings where Holmgren’s uniqueness  does not guarantee unique continuation. As a consequence of our analysis, unique continuation is nonetheless recovered for wave solutions observed on $(0,T) \times \omega$ with initial data supported in \( \mathscr{O} \).

The proof of 
our result combines a high-frequency observability estimate—based on the propagation of singularities—with a compactness-uniqueness argument that exploits the unique continuation properties of elliptic operators.

By duality, this observability result leads to 
controllability results for the wave equation, ensuring that the projection of the solution onto \( \mathscr{O} \) can be controlled by means of controls supported in $\omega$, with optimal spatial support.

We also present several extensions of the main result, including the case of boundary observations, as well as a characterization of the observable fraction of the energy of the initial data from partial measurements on $(0,T) \times \omega$. Applications to wave control are discussed accordingly.
\end{abstract} 
 \maketitle

\tableofcontents

\section{Introduction and first results}
\subsection{Problem formulation}

Let $\Omega$ be a bounded domain of $\R^{n}$, with boundary $\partial\Omega$ of class $\mathscr{C}^{\infty}$. We set 
$$ 
	\L = \R\times \Omega \quad \text{and } \quad \d\L = \R\times \d\Omega .
$$ 
We also introduce  $A=(a _{ij}(x))$,  a $n \times n$ matrix  of  $\mathscr{C}^{\infty}$ coefficients, symmetric,  uniformly definite positive on a neighborhood of $\overline{\Omega}$, and we denote $\Delta_{A}=\sum_{i,j=1}^{n}\partial_{x_{j}}(a _{ij}(x)\partial_{x_{i}} \cdot )$  the corresponding Laplacian operator. 

In addition, we will assume that the geodesics of $\overline{\Omega}$ with respect to the metric $(a^{ij}(x))_{ij}$ have no contacts of infinite order  with  $\d\Omega$. This is a standing assumption used  to define the global Melrose-Sj\"ostrand flow, see Section \ref{Geometry}.

We consider then the following wave equation:
\begin{equation}
\left\{ 
\begin{array}{ll}
	\partial _{t}^{2}u-\Delta_{A}u=0,\quad &\text{in } \,\,\L, 
	\\ 
	u(t,.)=0,\quad &\text{on } \,\,\d\L,
	\\
	(u(0,\cdot),\partial _{t}u(0,\cdot ))=(u_{0},u_{1}). & 
	\end{array}%
	\right.  
	\label{waveequation}
\end{equation}

We recall that, for $(u_0, u_1) \in H_{0}^{1}(\Omega)\times L^{2}(\Omega)$, the equation \eqref{waveequation} is well posed and admits a unique solution in the space $\mathscr{C}^{0}(\R,H_{0}^{1}(\Omega))\cap \mathscr{C}^{1}(\R,L^{2}(\Omega))$. It is also well posed in $L^{2}(\Omega)\times H^{-1}(\Omega)$ with the unique solution lying in $\mathscr{C}^{0}(\R,L^{2}(\Omega))\cap \mathscr{C}^{1}(\R,H^{-1}(\Omega))$.

The observation subdomain is denoted by $\omega$, a non-empty open subset of $\Omega$ where waves will be observed, and, for a given time-horizon $T>0$, we set $\omega_{T}=(0,T)\times\omega$ and $\L_{T}=(0,T)\times\Omega$. 

Thus, in the following, $\omega_T$ corresponds to the space-time subset on which $u$ is measured, observed and known, and our goal is to recover the initial data $(u_0, u_1)$ out of this partial measurement.

In other words, our aim is to analyse the inverse of the map
\begin{equation}
\label{Def-Obs-T}
	\begin{array}{lrll}
		\Obs_T: & L^2(\Omega) \times H^{-1}(\Omega) &\rightarrow & L^2(\omega_T)
		\\
		& (u_0, u_1) &\mapsto & u|_{\omega_T}, \text{ where $u$ solves \eqref{waveequation}.}
	\end{array}
\end{equation}
This is a classical question, motivated, in particular, by control theory, and commonly addressed in the context of the observability inequality 
\begin{equation}\label{Obs}
\Vert (u_{0},u_{1})\Vert_{L^2(\Omega) \times H^{-1}(\Omega)} \leq C \Vert u \Vert_ {L^{2}(\omega_{T})},  \end{equation}
which has been the focus of extensive research. This inequality refers to the possibility of inverting continuously the observation operator $\Obs_T$.

Traditionally it has been addressed for all possible solutions  with initial data in $L^2(\Omega) \times H^{-1}(\Omega)$ and it is  equivalent to the property of exact controllability of system \eqref{waveequation}
 in the space $H^1_0(\Omega)\times L^{2}(\Omega)$ with controls in $L^2(\omega_T)$, see \cite{Lions}. 
 
 Problem \eqref{Obs} is well understood, and there is  an almost necessary and sufficient condition for the validity of \eqref{Obs}, the so-called  Geometric Control Condition (in short \textit{GCC}, see \cite{Rauch-Taylor,B-L-R,BurqGerard}) which, roughly, asserts that all rays of Geometric Optics starting at time $t = 0$ from any point in $\overline{\Omega}$ meet the observation set $\omega_T$.
As we shall see below (see Section \ref{Geometry}) in more detail, these rays are the space-time projections of the generalized bicharacteristics  of Melrose-Sj\"ostrand,   \cite{MeSj}, for the operator $P_A = \partial_t^2 - \Delta_A$, which bounce on the boundary according to the Descartes-Snell law (see Section \ref{boundary-geometry}).
 
The GCC imposes a significant restriction on the class of observation domains $\omega$ for which the observability estimate \eqref{Obs} holds. In this work, we adopt a complementary viewpoint: rather than fixing the solution space and seeking suitable observation sets, we aim to consider all possible subdomains $\omega$, which is particularly relevant from an applied perspective due to practical constraints on the placement and availability of sensors or actuators. This naturally leads to the following question:
Can we characterize the subclass of solutions to \eqref{waveequation} for which the observability estimate \eqref{Obs} holds, given an arbitrary observation region $\omega$?

The main novelty of this paper lies in the sharp characterization of a class of solutions for which the observability estimate \eqref{Obs} holds without imposing any geometric conditions on the observation set $\omega$. Specifically, we show that \eqref{Obs} remains valid when the initial data of the solutions are supported in another subset \( \mathscr{O} \subset \Omega \), provided a suitable microlocal geometric condition is satisfied. This condition involves the time horizon $T$ and the pair \( (\omega, \mathscr{O}) \), and can be interpreted as a localized version of the Geometric Control Condition (GCC): it requires that all rays of geometric optics emanating from \( \mathscr{O} \) reach the observation region $ \omega$ within time $T$.
\subsection{Main results}\label{Sec-first results}

Our first main result is as follows:
\begin{theorem} \label{Theo-Obs-1}
	Given the domain $\Omega$, the observation subdomain $\omega \subset \Omega$, and the time-horizon $T>0$, let the subdomain $\mathscr{O}$  be such that
	every generalized bicharacteristic ray (see Definition \ref{def:gene-bichar} for their precise definition) starting from $\{t = 0\}\times \overline{\mathscr{O}}$  intersects the set $(0,T) \times \omega$.
	
	Then, there exists $C>0$ such that the solution of  \eqref{waveequation} satisfies the observability estimate  
	\begin{equation}
	\label{obs-theo-L2}
		\Vert (u_{0},u_{1})\Vert_{L^{2}(\Omega)\times H^{-1}(\Omega)} \leq C \Vert u\Vert_{L^{2}(\omega_{T} )},
			\end{equation}
		for  any initial data $(u_{0},u_{1}) \in L^{2}(\Omega)\times H^{-1}(\Omega)$ with support in $\overline{\mathscr{O}}$, i.e., satisfying
	\begin{equation}
	\label{Cond-Support-O}
		\supp(u_0, u_1) \subset	 \overline{\mathscr{O}}.
	\end{equation}

	Similarly, 
	 the observability estimate  
\begin{equation}\label{obs-theo}
\Vert(u_{0},u_{1})\Vert_{H^{1}_0(\Omega)\times L^{2}(\Omega)} \leq C \Vert \d_{t}u\Vert_{L^{2}(\omega_{T} )}
\end{equation}
holds for  any initial data $(u_{0},u_{1}) \in H_{0}^{1}(\Omega)\times L^{2}(\Omega)$ with support in $\overline{\mathscr{O}}$. 
\end{theorem}

\begin{rema}
This result holds for all subdomains $\omega$ and initial data with support in $\overline{\mathscr{O}}$, under the condition that the pair $(\omega, \mathscr{O})$ fulfills a mutual or joint microlocal relation, so that all rays emanating from $\overline{\mathscr{O}}$ reach $(0,T) \times \omega$. It is a natural extension of the classical one guaranteeing the observability of all solutions, namely the GCC, which corresponds to the case $\mathscr{O}=\Omega$ in our setting. Indeed, it is natural that, when the support of the initial data lies in $\overline{\mathscr{O}}$, its observation only depends on the dynamics of the rays emanating from $\overline{\mathscr{O}}$, independently of what other rays (the ones departing away from $\overline{\mathscr{O}}$) do. 

Our result extends the classical ones, allowing to consider all possible subdomains $\omega$, not only those fulfilling the highly demanding GCC, and identifying a class of observable initial data. This is particularly relevant in applications where the available observations are limited either because of the lack of accessibility to some regions of the domain where waves propagate or due to the lack of sufficient sensoring devices.	

   After this work has been submitted, the referee kindly pointed out to us the work \cite[Section 2]{Bardos-1993}, which deals with a similar problem in a particular setting in which $\mathscr{O}$ is a ball $B(r)$ strictly included in $\Omega$ and the metric is flat, under a slightly more demanding condition. Namely, the article \cite{Bardos-1993} requires that there is a neighborhood $B(r+\delta)$ ($\delta>0$) of the ball for which every generalized bicharacteristic ray passing through $[ 0,\delta ]\times B(r+\delta)$ intersects $(0,T) \times \omega$ (or $(0,T) \times \Gamma$ at a non-diffractive point if $\Gamma$ is a part of the boundary on which the observation is done, see Theorem \ref{Theo-Obs-2} for the counterpart of Theorem \ref{Theo-Obs-1} when the observation is performed on the boundary). Using partition of unity, we could derive from \cite{Bardos-1993} cases in which $\mathscr{O}$ is strictly included in $\Omega$, but it will not be sufficient to handle cases in which $\overline{\mathscr{O}}$ touches the boundary $\partial\Omega$ as in the example presented in Figure \ref{Fig-GCC-Interior} afterwards.
\end{rema}

\begin{rema} 
These results enter in the context of ``enlarged observability/controllability" introduced in \cite{Lions}, according to which when limiting the class of solutions under consideration the requirements for observability can be weakened. However, in \cite{Lions}, because of the use of the multiplier method, improvements were only achieved at the level of the needed observability time. The results of the present paper constitute a  more precise answer to these questions.
\end{rema}

A similar result holds when the observation is done along the boundary.  For this, we need the notion of nondiffractive points of the boundary that will be detailed in Definition \ref {def:nondiffractive}, see also \cite[Definition, pp.1037]{B-L-R} .

\begin{theorem} \label{Theo-Obs-2}
Let $\Gamma$ be a non-empty open subset of the boundary $\d\Omega$ and  $\mathscr{O}$ be a non-empty open set of  $\Omega$ such that the pair $(\Gamma, \mathscr{O}$) satisfies the following microlocal condition for some $T>0$: every generalized bicharacteristic ray starting from $\{t=0\}\times\overline{\mathscr{O}}$ intersects the set $(0,T)\times \Gamma$ at a nondiffractive point.
 
 Then there exists $C>0$ such that for  any initial data $(u_{0},u_{1}) \in H_{0}^{1}(\Omega)\times L^{2}(\Omega)$ supported in $\overline{\mathscr{O}}$, 
the solution of  \eqref{waveequation}  satisfies the observability estimate  
\begin{equation}\label{relaxed-obs-B}
\Vert (u_{0},u_{1})\Vert_{H_{0}^{1}(\Omega)\times L^{2}(\Omega)} \leq C \Vert \d_{n}u_{\vert\partial\Omega}\Vert_{L^{2}(\Gamma_{T} )} .
\end{equation}
\end{theorem}

\subsection{Examples}


\noindent \textbf{A $1$-d example.}
The microlocal assumption on the pair $(\omega, \mathscr{O})$ in Theorem \ref{Theo-Obs-1}, is, in general, weaker than the GCC since it only concerns the rays emanating from $\overline{ \mathscr{O}}$. It also requires a shorter observation  time. This is even the case in $1$-d.
	
\begin{figure}
\label{Fig-Ex-$1$-d}
\begin{tikzpicture}[scale=3]
    \draw[->] (-1.2, 0) -- (1.2, 0) node[below] {$x$};
    \draw[->] (0, 0) -- (0, 1.5) node[left] {$t$};
    \draw[blue,dashed] (-1,1) -- (1, 1) node[right] {$T = 1$}; 
    
    \draw[thick] (-1, 0) -- (-1, 1.5); 
    \draw[thick] (1, 0) -- (1, 1.5) ; 

    \draw[fill=blue, opacity=0.2] (-1, 0) rectangle (-0.75, 1.5);
    \draw[fill=blue, opacity=0.2] (0.75, 0) rectangle (1, 1.5);
    \node[above, blue] at (-0.875, 0) {$\omega$};
    \node[above, blue] at (0.875, 0) {$\omega$};

    \draw[fill=green, opacity=0.2] (-0.25, 0) rectangle (0.25, 1.5);
    \node[below, green!60!black] at (0, 0) {$\mathscr{O}$};

    \draw[red, thick, ->] (-0.25, 0) -- (1, 1.25) ;
    \node[below right] at (0.5,0.75) {{\color{red}$x_0 + t$}};
    \draw[red, thick, ->] (-0.25, 0) -- (-1, 0.75) ;
    \node[above] at (-0.5,0.4) {{\color{red}$x_0 - t$}};

    \draw[red, thick, dashed,->] (1, 1.25) -- (0.75, 1.5);
    \draw[red, thick, dashed,->] (-1, 0.75) -- (-0.5, 1.25);

    \node[below] at (-1, 0) {$-1$};
    \node[below] at (1, 0) {$1$};
    \node[below] at (-0.75, 0) {$-\frac{3}{4}$};
    \node[below] at (0.75, 0) {$\frac{3}{4}$};
    \node[below] at (-0.25, 0) {$-\frac{1}{4}$};
    \node[below] at (0.25, 0) {$\frac{1}{4}$};

\end{tikzpicture}
\caption{Illustration of the $1$-d example: $\Omega = (-1,1)$, $\omega = (-1,-3/4) \cup (3/4,1)$, $\mathscr{O} = (-1/4, 1/4)$. The critical time given by Theorem \ref{Theo-Obs-1} is $T_{0,crit} = 1$, while the times for unique continuation and the GCC coincide and are equal to $T_{UC} = 3/2$.}
\end{figure}
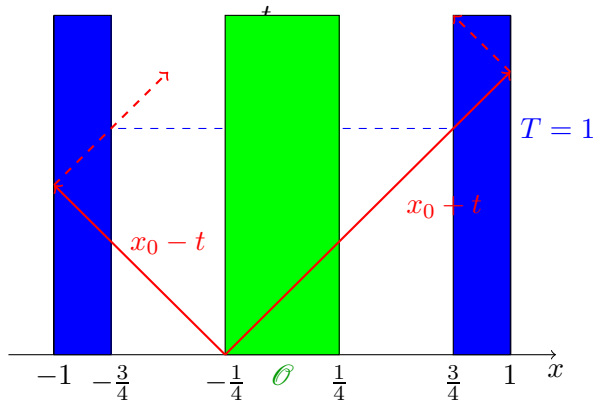

Indeed, consider the simple example of the $1$-d wave equation set on $\Omega = (-1,1)$, with control in $\omega = (-1, -3/4) \cup (3/4, 1)$, and initial data localized in $\mathscr{O} =(-1/4,1/4)$, as in Figure \ref{Fig-Ex-$1$-d}. In this case, the geodesic rays  are simply the characteristics $t \mapsto x_0 \pm t$, bouncing when meeting the boundaries $\{-1,1\}$.

By symmetry considerations, it is then easy to check that  Theorem \ref{Theo-Obs-1} holds for any $T > 1$. This minimal time corresponds to the arrival in $\omega$  of a characteristic starting from $x = -1/4$, propagating to the right. 

However, the classical, sharp condition for unique continuation in \eqref{Uniqueness-Condition} (in the following, $T_{UC}$ stands for the critical time for unique continuation; it is given by $2 \sup_{x \in \Omega} d(x, \omega)$ from the H\"ormander-Tataru-Robbiano-Zuily uniqueness theorem \cite{Tataru95, RobbianoZuily, Hormander-1997}) or the GCC require $T>T_{UC} = 3/2$ . This condition is indeed optimal  when aiming to observe all solutions since one can build initial data localized in $(-3/4, -3/4 + \varepsilon)$ ($\varepsilon >0$ small) leading to waves propagating towards the right at speed one, and vanishing in $\omega$ during the time interval $(0,3/2 - \varepsilon)$.

Therefore, the global observability estimates  \eqref{obs-theo-L2} or \eqref{obs-theo} do not hold for \emph{all} initial data for the intermediate times $1 < T < 3/2$, but, according to our result, they do hold for initial data localized in $\overline{\mathscr{O}}$. 

\medskip

\noindent \textbf{Multi-d examples.}
	We present below two additional multi-d examples:
\begin{enumerate}[leftmargin=0pt]
	\item Let us consider the case illustrated in Figure \ref{Fig-Ex-Ball-multid}, in which $\Omega$ is the unit ball, $\omega$ is the $\varepsilon$-neighborhood of its boundary and $\mathscr{O}$ is the interior ball centred at the origin, of radius $\alpha$, with $0 < \alpha < 1 - \varepsilon$. Then the critical time given by Theorem \ref{Theo-Obs-1} is  $T_{0,crit} =   1+\alpha - \varepsilon$. But, in this case, the time for unique continuation and GCC is larger, namely,  $2(1 -  \varepsilon)$.
	\begin{figure}
	\label{Fig-Ex-Ball-multid}
	\begin{tikzpicture}[scale=2]
    \draw[thick] (0, 0) circle (1);
    \node at (0, 0) {$\Omega$};

    \draw[thick, blue] (0, 0) circle (0.8);
    \fill[blue, opacity=0.2] (0, 0) circle (1);
    \fill[white] (0, 0) circle (0.8);
    \node[right, blue] at (1, -0.3) {$\omega$};

    \draw[thick, green!70!black] (0, 0) circle (0.3);
    \fill[green, opacity=0.3] (0, 0) circle (0.3);
    \node[left, green!70!black] at (-0.2, 0.3) {$\mathscr{O}$};

    \draw[<->, thick] (0, 0) -- (0,1) node[midway, right] {$1$};
    \draw[<->, thick, blue] (0.8, 0) -- (1, 0) node[midway, above ] {$ \varepsilon$};
    \draw[<->, thick, green!70!black] (0, 0) -- (0.3, 0) node[midway, below] {$\alpha$};

	\draw[red, thick, <->] (-0.15, -0.26) -- (0.5,0.86) node[right]{A critical characteristic} ;
    
\end{tikzpicture}
\caption{llustration of the multi-d example (Item 1): $\Omega = B(1)$, $\omega = B(1)\setminus B(1- \varepsilon)$, $\mathscr{O} = B(\alpha)$. The critical time given by Theorem \ref{Theo-Obs-1} is $T_{0,crit} = 1+ \alpha - \varepsilon$, while the times for unique continuation and GCC coincide: $T_{UC} = 2(1- \varepsilon)$.}
\end{figure}
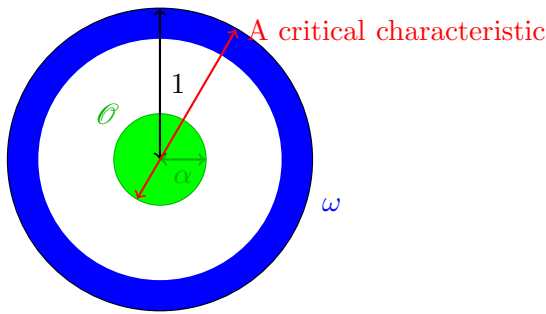

\item Another example, illustrated in Figure \ref{Fig-GCC-Interior}, still when $\Omega$ is the unit ball of $\R^2$, corresponds to an observation  subdomain $\omega$ which is an $\varepsilon$ neighborhood of one third of the boundary of $\Omega$, with angles in $(-\pi/3,  \pi/3)$, and $\mathscr{O}$ being the vertical strip $\{ x = (x_1,x_2) \in \Omega \hbox{ with } x_1 < \cos(\alpha) \} $ for some $\alpha \in (\pi, 4 \pi/3)$. The critical time $T$ given by Theorem \ref{Theo-Obs-1} is finite but GCC fails, whatever the time-horizon is, due to the  vertical diameter, corresponding to a geodesic ray which bounces back and forth endless, without ever entering the observation set $\omega$.

The longest geodesic that starts from $\mathscr{O}$ and stays away from $\omega$ is the one that starts from $(\cos(\alpha), \sin(\alpha))$ and goes to $(\cos(\pi/3), \sin(\pi/3))$ (or rather an $\varepsilon$-neighborhood of it): it is not difficult to check that this geodesic has length $4 k_c \sin( (\alpha - \pi/3)/2)$ where $k_c$ is the first integer such that $2 k_c (\pi - (\alpha- \pi/3) ) > 5 \pi /3 - \alpha$.

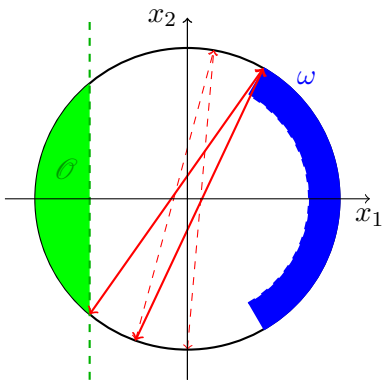
\begin{figure}[h!]
\label{Fig-GCC-Interior}

\begin{tikzpicture}[scale=2]
    \draw[thick] (0, 0) circle (1);

    \fill[blue, opacity=0.2, domain=-60:60, samples=100] 
        plot ({cos(\x)}, {sin(\x)}) 
        -- plot ({0.8*cos(60)}, {0.8*sin(60)}) arc (60:-60:0.8) -- cycle;
    \draw[blue, thick, domain=-60:60, samples=100] 
        plot ({cos(\x)}, {sin(\x)}) 
        node[blue, right] at (0.65, 0.8) {$\omega$};

    \draw[blue, thick, dashed, domain=-60:60, samples=100] 
        plot ({0.8*cos(\x)}, {0.8*sin(\x)});

    \fill[green, opacity=0.3] 
    ({cos(130)},{sin(130)}) -- ({cos(230)}, {sin(230)}) arc[start angle=230, end angle=130, radius=1] -- cycle;
    \node[green!60!black] at (-0.8, 0.2) {$\mathscr{O}$};

    \draw[thick, dashed, green!70!black] ({cos(130)}, -1.2) -- ({cos(130)}, 1.2);

    \draw[red, thick, <->] ({cos(230)}, {sin(230)}) -- ({cos(60)}, {sin(60)});
    \draw[red, thick, <->] ({cos(250)}, {sin(250)}) -- ({cos(60)}, {sin(60)});
    \draw[red, dashed, <->] ({cos(250)}, {sin(250)}) -- ({cos(80)}, {sin(80)});
    \draw[red, dashed, <->] ({cos(270)}, {sin(270)}) -- ({cos(80)}, {sin(80)});

    \draw[->] (-1.2, 0) -- (1.2, 0) node[below] {$x_1$};
    \draw[->] (0, -1.2) -- (0, 1.2) node[left] {$x_2$};

    \draw[->] (0, -1.2) -- (0, 1.2); 
\end{tikzpicture}
   \caption{Illustration of the multi-d example (Item 2),  with $\alpha = 230^\circ$: In this case, it is clear that the critical time $T$ given by Theorem \ref{Theo-Obs-1} is finite, while the Geometric Control Condition is never satisfied. In red, we have plotted the first few reflections of the longest geodesic starting from $\mathcal{O}$ and staying outside $\omega$.}
\end{figure}

\end{enumerate} 

\begin{rema}{The examples above can be easily adapted to the boundary-observation setting by ``smashing'' the set $\omega$ to the boundary of $\Omega$.}
\end{rema}

\subsection{Some relevant consequences}
As we  will see in Section \ref{Sec-Control-without-GCC} below,  the observability results in  Theorems \ref{Theo-Obs-1} and \ref{Theo-Obs-2}, have their dual controllability counterparts: One can exactly control the projection of wave solutions over $\mathscr{O}$ at time $T$, by means of controls in $L^2(\omega_{T})$, see Theorem \ref{Thm-Obs-local}. 

Such result lays in between the classical properties of approximate and exact controllability, since it assures the exact control over the projection onto $\mathscr{O}$, but without providing any information on what happens outside $\mathscr{O}$. 

Note that, if we further assume that the time horizon $T$ is sufficiently large to guarantee the injectivity of the operator $\Obs_T$ in \eqref{Def-Obs-T} (in fact, this corresponds to the condition $T > 2 \sup d(x,\omega)$, see \eqref{Uniqueness-Condition} afterwards), we can find controls that simultaneously assure the exact control of the projection over $\mathscr{O}$ and the approximate controllability property everywhere in the domain $\Omega$, see Remark \ref{Remark-Exact-Control+GlobalApprox}.

The problem of  controllability in the absence of GCC has been the object of extensive study, see, for instance,  \cite{Russell-1971-Wave, Lebeau92, Robbiano95, Laurent-Leautaud}. In those papers, one aims to quantify the property of approximate controllability, by identifying subspaces of initial data that can actually be controlled. These spaces are typically very small, imposing suitable analyticity restrictions. 
Surprisingly, analyzing what can be controlled for solutions of the wave equation in the usual energy space in the absence of GCC has not received much attention, and, besides our work, we are only aware of the work \cite{Bardos-1993} focusing on a very particular case.


Our result, valid for arbitrary observation sets $\omega$, encompasses as a particular case the classical setting in which \( \mathscr{O} = \Omega \), that is, when observability is required for all solutions of the wave equation. In this case, our condition naturally reduces to the well-known Geometric Control Condition (GCC).

Our main result, as we shall see, is even sharper, since it allows to identify the observable microlocal projections of solutions for any subdomain $\omega$, something that might be also relevant in applications. 

By duality, these observation type properties yield control results for the wave equation in which the control objective is to control the projection of the solutions over $\mathscr{O}$, or its microlocal projection up to some smoother error terms (see Section \ref{Sec-Control-without-GCC} for the precise statement). 

Throughout the paper, we also present several illustrative examples involving different domains $ \Omega$, observation subsets $\omega$, and time horizons $T$, which do not satisfy the classical GCC. Nonetheless, our results apply in those cases, as our goal is to discuss all the information on the solutions of the wave equation that can be derived from a given observation set $(0,T)\times \omega$.

 Our analysis has also important consequences in what concerns the classical property of unique continuation. Indeed, the injectivity of the operator $\Obs_T$ is equivalent to the unique continuation property 
 \begin{equation}
 	\label{Unique-Continuation-Prop}
	\text{For $u$ solution of \eqref{waveequation}}, \,
	u|_{\omega_T} = 0 \Rightarrow (u_0, u_1) \equiv (0,0).
 \end{equation}
Such property is known to hold for all possible solutions of the wave equation when $A$ is analytic, thanks to Holmgren's uniqueness theorem, \cite{John-1949}, or for smooth time-independent coefficients, by \cite{Tataru95, RobbianoZuily, Hormander-1997}, provided the time $T$ satisfies 
\begin{equation}
	\label{Uniqueness-Condition}
	T > 2 \sup_{x \in \Omega} d(x, \omega),
\end{equation}
where $d(\cdot, \omega)$ stands for the geodesic distance to $\omega$. Accordingly,  unique continuation holds for all non-empty open subset $\omega$ provided $T$ is large enough as in \eqref{Uniqueness-Condition}. In this setting, a quantitative logarithmic stability estimate was recently proved in \cite{Laurent-Leautaud}.  Thus, under the only condition \eqref{Uniqueness-Condition}, the operator $\Obs_T$ is one-to-one, but in general, its inverse is ill-posed and it is not a bounded operator, unless  the additional \textit{GCC} is satisfied. 

As an interesting corollary of our observability result, we derive unique continuation properties for specific classes of solutions — such as those with initial data supported in \( \mathscr{O} \) — even in settings where the classical condition \eqref{Uniqueness-Condition} fails. Consequently, these results apply in situations where existing uniqueness theorems of Holmgren do not suffice.

Our contribution may also be seen as a complement or alternative to the results developed in \cite{Lebeau92} and the subsequent literature. In that context, for general subdomains $\omega$, and under the sole assumption of the unique continuation condition \eqref{Uniqueness-Condition}—in the absence of the GCC—observability for general solutions was established in a weaker, generalized framework, where the observability constant depends, roughly, exponentially on the frequency of the solutions.

In contrast, our approach aims to recover the classical observability inequality in the natural energy spaces, even without the GCC, by identifying—microlocally—specific classes of initial data for which the inequality still holds. Rather than relaxing the inequality to include all solutions at the cost of weakening the estimate, we preserve its sharp form within a suitably restricted solution space. Throughout the paper, we clarify the relationship between both approaches and highlight the connections between their respective results.

\subsection{Methodology of proof}\label{methodology}
The proof strategy for Theorems \ref{Theo-Obs-1} and \ref{Theo-Obs-2} relies on microlocal analysis techniques, as is customary in the study of wave propagation phenomena. More precisely, it combines the following key ingredients:
\begin{itemize}
\item The microlocal geometric property satisfied by the pair $(\omega, \mathscr{O})$ (or its boundary counterpart $(\Gamma, \mathscr{O})$) allows to propagate the energy of solutions at high-frequencies from $(0,T) \times \omega$ towards $\{ t = 0 \} \times \mathscr{O}$, and this leads to a weak version of the observability estimate with a compact remainder term. 
\item Removing the compact remainder requires  a unique continuation property. This property needs to be proved in an ad-hoc manner since the assumptions made on the time-horizon do not assure that \eqref{Uniqueness-Condition} is fulfilled.
\item This is achieved by means of an added compactness-uniqueness argument, which reduces the unique continuation property to an elliptic context in which it holds by classical Carleman inequalities.
\end{itemize}

\subsection{Outline.} \label{outline}
Section \ref{Proofs 1-2} is devoted to the proofs of the main results—Theorems \ref{Theo-Obs-1} and \ref{Theo-Obs-2}—following the methodology outlined in Section \ref{methodology}. In Section \ref{Geometry}, we recall the Melrose–Sj\"ostrand cotangent bundle framework and introduce several technical microlocal analysis tools that are used throughout the paper. Section \ref{sec-further results} presents extensions of the observability results, where the assumptions on the support of the initial data are relaxed. In Section \ref{Sec-Control-without-GCC}, we establish the control counterparts of the observability results. The article concludes with a final section discussing open problems and future perspectives.
\bigskip

\noindent{\it Acknowledgments.} The authors thank the referee for having pointed out the article \cite{Bardos-1993}.

\section{Proofs of the main results}\label{Proofs 1-2}
We essentially focus on the proof of Theorem \ref{Theo-Obs-1}, i.e., of the estimate \eqref{obs-theo}, since the proof of estimate \eqref{obs-theo-L2} is similar, and is thus left to the reader, with the indication of the
additional steps needed, see Section \ref{Subsec-Proof-Theo-Obs-2} afterwards. The general strategy of the proof follows the program described in Section \ref{methodology}.

\subsection{Proof of Theorem \ref{Theo-Obs-1} }
We start with the following lemma that describes the propagation of regularity for solutions to system  \eqref{waveequation}. The proof of this lemma requires the use of microlocal tools and is a consequence of \cite{MeSj,Hormander-III}, and we refer the unfamiliar reader to Section \ref{Geometry} for a precise description of these notions. In particular, the proof of Lemma \ref{regularity} can be found in Section \ref{Subsec-Proof-Lemma-reg}.

\begin{lemma}\label{regularity}
Under the assumptions of Theorem \ref{Theo-Obs-1}, consider  a solution $u$ to system \eqref{waveequation} with initial data $(u_{0},u_{1}) \in H^{-1}(\Omega)\times (H^{2}\cap H^1_0)'(\Omega)$ supported in $\overline{\mathscr{O}}$,  and  satisfying $u \in L^{2}(\omega_{T} )$. 
Then  $(u_{0},u_{1}) \in L^{2}(\Omega)\times H^{-1}(\Omega)$ and $u \in \mathscr{C}^0(\R, L^2(\Omega)) \cap \mathscr{C}^1(\R, H^{-1}(\Omega))$.

Similarly, if  $(u_{0},u_{1}) \in L^{2}(\Omega)\times H^{-1}(\Omega)$ is supported in $\overline{\mathscr{O}}$,  and  the solution $u$ satisfies $\d_{t}u \in L^{2}(\omega_{T} )$, one has  $(u_{0},u_{1}) \in H_{0}^{1}(\Omega)\times L^{2}(\Omega)$  and $u \in \mathscr{C}^0(\R, H^1_0(\Omega)) \cap \mathscr{C}^1(\R, L^2(\Omega))$.
\end{lemma}

We then deduce the following corollary, which fulfills the first step of the proof, providing a first observability estimate with a compact remainder.

\begin{corollary}\label{relaxed-obs}
Under assumptions of Theorem \ref{Theo-Obs-1},  there exists $C>0$ such that 
the solution of  \eqref{waveequation}  satisfies 
\begin{equation}\label{R-obs-L2}
\Vert (u_{0},u_{1})\Vert_{L^{2}(\Omega)\times H^{-1}(\Omega)} \leq C \Vert u\Vert_{L^{2}(\omega_{T} )} +C \Vert(u_{0},u_{1})\Vert_{H^{-1}(\Omega)\times (H^{2}\cap H^1_0)'(\Omega)},
\end{equation}
for  any initial data $(u_{0},u_{1}) \in L^{2}(\Omega)\times H^{-1}(\Omega)$ supported in $\overline{\mathscr{O}}$ (i.e., satisfying \eqref{Cond-Support-O}).

Similarly, 
\begin{equation}\label{R-obs-L1}
\Vert(u_{0},u_{1})\Vert_{H^{1}_0(\Omega)\times L^{2}(\Omega)} \leq C \Vert \d_{t}u\Vert_{L^{2}(\omega_{T} )} +C \Vert(u_{0},u_{1})\Vert_{L^{2}(\Omega)\times H^{-1}(\Omega)},
\end{equation}
holds for  any initial data $(u_{0},u_{1}) \in H_{0}^{1}(\Omega)\times L^{2}(\Omega)$  supported in $\overline{\mathscr{O}}$ (i.e., satisfying \eqref{Cond-Support-O}).
\end{corollary}

\begin{proof}
 We focus on the proof of  \eqref{R-obs-L1}, \eqref{R-obs-L2} being similar.
Consider  the following Hilbert space 
\begin{equation*}
	E=\Big\{ (u_{0},u_{1}) \in L^{2}(\Omega)\times H^{-1}(\Omega), \, \text{supp}(u_{0},u_{1}) \subset \overline{\mathscr{O}}, \,\, \text{ and } \d_{t}u \in L^{2}(\omega_{T})\Big\}
\end{equation*}
equipped with the norm 
$$\Vert (u_{0},u_{1})\Vert_{E}^{2}=\Vert (u_{0},u_{1})\Vert_{L^{2}(\Omega)\times H^{-1}(\Omega)}^{2} + \Vert \d_{t}u\Vert_{L^{2}(\omega_{T})}^{2},$$ 
and the  energy space $F = H_{0}^{1}(\Omega)\times L^{2}(\Omega)$
equipped with its natural norm. 

Thanks to Lemma \ref{regularity},  the identity map 
\begin{equation*}
	E \longrightarrow F = H_{0}^{1}(\Omega)\times L^{2}(\Omega) ,  \quad (u_{0},u_{1}) \mapsto (u_{0},u_{1})
\end{equation*}
is well defined. Consequently,  the closed graph theorem  yields its continuity and estimate \eqref{R-obs-L1}. 

	The proof of estimate \eqref{R-obs-L2} is similar, except that it relies on the propagation of the $L^{2}$-wave front set along the bicharacteristic flow, that is the first item of Lemma \ref{regularity}.
\end{proof}

As a second step in the proof, and as a consequence of the previous estimate, the following unique continuation property holds:

\begin{lemma} \label{UC}
Under assumptions  of Theorem \ref{Theo-Obs-1},  any solution $u$ of system \eqref{waveequation} with initial data in $H^{-1}(\Omega)\times (H^{2}\cap H^1_0(\Omega))'$ with support in $\overline{\mathscr{O}}$,  and satisfying $u= 0$ in $ \omega_{T} $,  vanishes everywhere, i.e., $u \equiv 0$.

Similarly, any solution $u$ of system \eqref{waveequation} with initial data in $L^2(\Omega)\times H^{-1}(\Omega)$ supported in $\overline{\mathscr{O}}$,  and satisfying $\d_{t}u= 0$ in $ \omega_{T} $,  vanishes everywhere, i.e., $u \equiv 0$.
\end{lemma} 

\begin{rema} At this point it is worth noticing that this uniqueness result is not standard since it only applies to the solutions with initial data supported in $\overline{\mathscr{O}}$. It is not a consequence of Holmgren's uniqueness theorem nor any of its generalisations, but rather a corollary of Corollary \ref{relaxed-obs}, which establishes a relaxed version of the observability inequalities we aim, with an added compact additive remainder term.
\end{rema}
\begin{proof} 
	Similarly as in the proof of Corollary \ref{relaxed-obs}, we focus on finite energy solutions, the case of weaker solutions being similar.
	
Our goal is to prove that the closed linear subspace of  $L^{2}(\Omega) \times H^{-1}(\Omega)$ defined by
$$
\mathcal{N} = \Big\{  (u_{0},u_{1}) \in L^2(\Omega)\times H^{-1}(\Omega), \, \text{Supp}(u_{0},u_{1}) \subset \overline{\mathscr{O}} ,\, \, \d_{t}u_{\vert\omega_{T}} =0 \Big\} 
$$
is reduced to $\mathcal{N}=\{(0,0)\}$.

Thanks to the estimate \eqref{R-obs-L1} in Corollary \ref{relaxed-obs},   it is clear that $\mathcal{N} \subset H_{0}^{1}(\Omega) \times L^2(\Omega)$ and \begin{equation*}
\Vert(u_{0},u_{1})\Vert_{H^{1}_0(\Omega) \times L^{2}(\Omega)} \leq C \Vert(u_{0},u_{1})\Vert_{L^{2}(\Omega)\times H^{-1}(\Omega)}
\end{equation*}
for every $(u_{0},u_{1}) \in \mathcal{N}$.

Using then the  compact embedding $H^{1}_{0}(\Omega) \hookrightarrow L^{2}(\Omega)$ \footnote{In fact, since we are considering data which are supported in $\overline{\mathscr{O}}$, we can also consider the compact embedding of $H^{1}_{0}(\mathscr{O}) \hookrightarrow L^{2}(\mathscr{O})$. This remark allows to deduce the same results even in cases in which $\Omega$ is unbounded, provided $\mathscr{O}$ is bounded.} we deduce that $\mathcal{N}$ has a finite dimension. In addition,   the matrix operator ${\mathcal A} = \begin{pmatrix} 0 &1\\ \Delta_{A} & 0\end{pmatrix}$ defines a linear bounded operator in $\mathcal{N}$. This is so since the wave equation, written on the column vector unknown $U=\begin{pmatrix} u\\ u_t\end{pmatrix}$ takes  the form $U_t= {\mathcal A}U$. Therefore,  ${\mathcal A}$ operates continuously in $\mathcal{N}$ and corresponds to the application of the time derivative on the solutions of the wave equation \eqref{waveequation}, transferring the initial data $(u_{0},u_{1})$ into $(u_{1},\Delta_Au_{0}).$ 

If $\mathcal{N}$ were non-trivial, the operator ${\mathcal A}$ would have an eigenvalue corresponding to an eigenvector in $\mathcal{N}$. But, as we shall see, this is impossible, concluding, by contradiction, that $\mathcal{N}=\{(0,0)\}$. 

Indeed, if $\lambda \in \C$ is an eigenvalue and $(u_{0}, u_1)\ne (0, 0)$  in $H_{0}^{1}(\Omega) \times L^2(\Omega)$ is an eigenvector  of $\mathcal{A}$ in $\mathcal{N}$, we have,
$$ 
\Delta_{A} u_{0} -\lambda^{2}u_{0}=0  \text{ in } \Omega, \qquad u_0 = 0 \text{ on }  \partial \Omega, \qquad \hbox{ and } \qquad \lambda u_{0} =0  \text{ in }  \omega.
$$
It is easy to check that this cannot happen unless $u_0 \equiv 0$, which also implies that $u_1\equiv 0$. Indeed, if $\lambda =0$, given that $\Delta_{A} u_{0} =0$ and $u_0$ vanishes on the boundary, we conclude that $u_{0} \equiv 0$.   On the other hand, when $\lambda \ne 0$, $u_{0} $ vanishes in $\omega$ and by elliptic unique continuation applied to the equation $ 
\Delta_{A} u_{0} -\lambda^{2}u_{0}=0$ we deduce that $u_{0}  \equiv 0$ everywhere.  

This concludes the proof of Lemma \ref{UC}.
\end{proof} 

We are now in conditions to conclude the proof of Theorem \ref{Theo-Obs-1}, i.e., of estimate \eqref{obs-theo}.

We use a contradiction argument and we assume that estimate \eqref{obs-theo} is false. Consider  a sequence of initial data $(u_{0,k},u_{1,k}) \in H_{0}^{1}(\Omega)\times L^{2}(\Omega)$ with support in $\overline{\mathscr{O}}$,  and $(u_{k})$ the corresponding solutions, with
\begin{equation}\label{contradiction1}
	\Vert (u_{0,k},u_{1,k})\Vert_{H^{1}_0(\Omega) \times L^{2}(\Omega)} = 1,
	  \quad\text{and} \quad 
	  \lim_{k \to \infty} \Vert \d_{t}u_{k}\Vert_{L^{2}(\omega_{T})} = 0 .
\end{equation}

The sequence $(u_{k})$ is bounded in the energy space $\mathscr{C}^{0}([0,T]; H^{1}_0(\Omega)) \cap \mathscr{C}^{1}([0,T]; L^{2}(\Omega))$. Thus, after extracting a subsequence,  we may assume that it converges weakly in $H^1(\L_{T})$ to another solution $u \in H^{1}(\L_{T})$ of \eqref{waveequation}, corresponding to an initial datum $(u_{0},u_{1}) \in H^1_0(\Omega) \times L^2(\Omega)$ with support in $\overline{\mathscr{O}}$, weak limit of $(u_{0,k},u_{1,k})$ in the energy space.

Passing  to the limit $k\rightarrow \infty$  in \eqref{contradiction1}, we obtain
\begin{equation}\label{weaklimit}
 \d_{t}u_{\vert \omega_{T}}= 0 .
\end{equation}
Then, the unique continuation result of Lemma \ref{UC}  assures that  $u\equiv 0$. This implies that $(u_{0,k},u_{1,k})$ strongly converges to $(0,0)$ in $L^{2}(\Omega)\times H^{-1}(\Omega)$.

On the other hand, in view of the  relaxed observability estimate \eqref{R-obs-L1} applied to $u_{k}$, we have 
\begin{equation*}
1 = \Vert (u_{0,k},u_{1,k})\Vert_{H^{1}_0(\Omega)\times L^{2}(\Omega)}  
\leq  C \Vert \d_{t}u_{k}\Vert_{L^{2}(\omega_{T})} + C \Vert (u_{0,k},u_{1,k})\Vert_{L^{2}(\Omega)\times H^{-1}(\Omega)}.
\end{equation*}
But the right hand-side tends to $0$ as $k\to \infty$ thanks to \eqref{contradiction1} and the fact that $(u_{0,k},u_{1,k}) \to (0,0)$ in $L^{2}(\Omega)\times H^{-1}(\Omega)$. \footnote{Note that, here again, the same argument applies if $\Omega$ is unbounded but $\mathscr{O}$ is bounded, since  the compact embedding of $H^{1}_0(\mathscr{O}) \times L^{2}(\mathscr{O})$ into $L^{2}(\mathscr{O}) \times H^{-1}(\mathscr{O})$ can be employed.}  This yields  a contradiction. 

\subsection{Proof of Theorem \ref{Theo-Obs-2}}\label{Subsec-Proof-Theo-Obs-2}
	To conclude this section, we outline the proof of Theorem \ref{Theo-Obs-2}, which closely follows that of Theorem \ref{Theo-Obs-1}. The only difference lies in the replacement of Lemma \ref{regularity} with the following lemma, whose proof is sketched in Remark \ref{Rema-Lem-regularity-boundary} in Section \ref{Subsec-Proof-Lemma-reg}:
	
	\begin{lemma}\label{regularity-boundary}
Under the assumptions of Theorem \ref{Theo-Obs-2}, consider  a solution $u$ to system \eqref{waveequation} with initial data $(u_{0},u_{1}) \in L^{2}(\Omega)\times H^{-1}(\Omega)$ supported in $\overline{\mathscr{O}}$,  and  satisfying $\partial_n u \in L^{2}(\Gamma_{T} )$. 
Then  $(u_{0},u_{1}) \in H_{0}^{1}(\Omega)\times L^{2}(\Omega)$ and $u \in \mathscr{C}^0(\R, H^1_0(\Omega)) \cap \mathscr{C}^1(\R, L^2(\Omega))$. 
\end{lemma}


\section{Some geometric facts, operators and wave fronts}\label{Geometry}
In this section, we analyse the geometry of the domain $\Omega$ near the boundary and we provide the microlocal material used in this paper. More precisely, we  present the generalized bicharacteristic flow of Melrose-Sj\"ostrand,  the notion of wave front set up to the boundary and the theorem of propagation of singularities. All these notions are borrowed to Melrose-Sj\"ostrand  \cite{MeSj} and H\"ormander \cite{Hormander-III}.

Recall that the compressed cotangent bundle of Melrose-Sj\"ostrand  is given by
\begin{equation*}
T^*_b\L=T^{\ast }\L\cup T^{\ast }\partial \L,
\end{equation*}
 and that we have a natural projection
\begin{equation}
\label{projection1}
\pi\ :\ T^{*}\R^{n+1}\mid _{\overline \Omega} \,\rightarrow T^{*}_{b}\L.
\end{equation}
We equip $T^{\ast }_b\L$ with the induced topology. 

Given the matrix $A(x)=(a_{ij}(x))$ and $\xi \in \R^{n}$, we set $\vert\xi\vert_{x}^{2}=^t\xi A(x)\xi$. We also denote by $p_A(t,x;\tau,\xi)= -\tau^2+\vert\xi\vert_{x}^{2}$, the principal symbol of  $P_{A}=\d_{t}^{2}-\sum_{i,j=1}^{n}\partial_{x_{j}}(a _{ij}(x)\partial_{x_{i}} \cdot )$. Finally, we define the characteristic set 
 \begin{equation*}
 \text{Char}(P_{A}) =  \{(t,x;\tau,\xi) \in T^{*}\R^{n+1}\backslash 0, \, p_{A}(t,x;\tau, \xi) =  0\},
\end{equation*}
and  $\Sigma_{A} = \pi ( \text{Char}(P_{A}))$. 

\subsection{Local geodesic coordinates }\label{geo}
Near a point  $m_{0}$ of the boundary $\d\Omega$, taking advantage of the regularity of $\Omega$ which is of class $\mathscr{C}^{\infty}$, we can define a system of geodesic coordinates $x=(x_{1},x_{2},....,x_{n}) \mapsto y=(y_{1},y_{2},....,y_{n})$ such that we have locally
 \begin{equation}\label{Geodesic-1}
\Omega= \{ (y_{1},y_{2},....,y_{n}), \,\, y_{n}>0\}, \quad \d\Omega = \{ (y_{1},y_{2},....,y_{n-1},0)\}=\{ (y',0)\} ,
 \end{equation}
and the corresponding wave operator is given by
\begin{equation*}
P_{A} = \partial_{t}^{2}-\Big(\partial_{y_{n}}^{2}+\sum_{1\leq i,j \leq n-1}\partial_{y_{j}}\left(b_{ij}(y)\partial_{y_{i}} \, \right)\Big) +M_{0}(y)\partial_{y_{n}} + M_{1}(y,\partial_{y'}) .
\end{equation*}
Here, the matrix $(b_{ij}(y))_{ij}$ is of class $\mathscr{C}^{\infty}$, symmetric,  uniformly definite positive on a neighborhood of $m_{0}$, $M_{0}(y)$ is a real valued function of class $\mathscr{C}^{\infty}$, and $ M_{1}(y,\partial_{y'})$ is a tangential differential operator  of order $1$ with $\mathscr{C}^{\infty}$ coefficients.

In the sequel, for convenience, we will use the same notation $(t,x)=(t,x',x_{n})$ to denote $(t,y',y_{n})$, and we shall write
 \begin{equation}\label{Geodesic-2}
P_{A} = -\partial_{n}^{2}-R(x_{n},x',\d_{x',t}) + M_{0}(x)\partial_{n} + M_{1}(x,\partial_{x'}).
 \end{equation}
Notice that, in this system of coordinates,  the principal symbol of the wave operator $P_{A}$ is given by 
\begin{equation*}
p_{A} = \xi_{n}^{2}-r(x,\tau,\xi') = \xi_{n}^{2} - \Big(\tau^{2}-\sum_{1\leq i,j \leq n-1}a_{ij}(x)\xi_{i}\xi_{j}\Big).
\end{equation*}
We shall set $r_{0}(x',\tau,\xi')=r(x',0, \tau,\xi')$.

\subsection{Generalized bicharacteristic rays}\label{boundary-geometry}

First, let us recall that the hamiltonian field associated to $p_A$ is given by 
\begin{equation*}
H_{p_A} = -2\tau \partial_t +2 ^t\xi A(x)\nabla_x - \sum_{k=1}^{n} \, ^t \xi(\partial_{x_k}A(x))\xi\partial_{\xi_k} .
\end{equation*}

We also recall the following partition of $T^{\ast}(\d\L)$ into elliptic, hyperbolic and glancing sets:
 \begin{equation}\label{Boundary}
 \#\Big\{\pi^{-1}(\rho) \cap  \text{Char}(P_A)\Big\} = 
\left\{ 
\begin{array}{c}
0 \quad if  \quad \rho \in \E
 \\ 
1 \quad if  \quad \rho \in \G
\\ 
2 \quad if  \quad \rho \in \H.
\end{array}%
\right.  
\end{equation}
For the sake of simplicity, we develop the rest of this section in a system of local geodesic coordinates as introduced in Section \ref{geo}.  Therefore  we have locally
$$
\E = \{r_{0}<0\}, \quad\quad\H = \{r_{0}> 0\}, \quad\quad\G = \{r_{0}=0\} .
$$
In addition, using the projection $\pi$, one can identify the glancing set $\G$ with a subset of $T^{\ast }\R^{n+1}$. 

Following \cite{Hormander-III} and 
\cite{B-L-R}, we have the precise description of the glancing set $\G$.

\begin{defi}\label{def:nondiffractive}
Let $\rho$ be a point of $T^{\ast }\partial \L\backslash 0$. We say that

\begin{enumerate}

\item $\rho $ is diffractive  if $\rho \in \G$ and $H_{p_A}^2(x_{n})(\rho) > 0$.

This means that the free bicharacteristic ray $\gamma$ issued from $\rho$ belongs to $T^{\ast }\L$ in a neighborhood of $0$, except at $s = 0$, i.e., there exists $\varepsilon >0$ such that $(\exp sH_{p_A})\widetilde{\rho } \in T^{\ast }\L$ for $0<\vert s\vert <\varepsilon$, with $\widetilde{\rho }=\pi^{-1}(\rho)$.

\item $\rho$ is nondiffractive if  a) $\rho \in \H$  or b) $\rho \in \G$  and the free bicharacteristic ray $(\exp sH_{p_A})\widetilde{\rho }$ passes over the complement of $\L$ for arbitrarily small values of $s$ .


\end{enumerate}
\end{defi}

We shall denote by $\G_d$ the set of diffractive points . Notice that in local geodesic coordinates, the set $\G_{d}$  is given by 
\begin{equation}\label{diffractive}
\G_{d}= \{\xi_{n}= r_{0}=0, \, \partial_{n}r_{\vert x_{n}=0}>0\} .
\end{equation}

\begin{defi}\label{def:gene-bichar}
A generalized bicharacteristic ray is a continuous map
$$
\R \supset I \setminus B\ni s \mapsto \gamma (s) \in T^{\ast }\L \cup \G \subset T^{\ast }\R^{n+1}
$$
where $I$ is an interval of $\R$, B is a set of isolated points, for every $s  \in I \setminus B$, $\gamma (s) \in \Sigma_A$ and $\gamma$ is differentiable as a map with values in $T^{\ast }\R^{n+1}$, and 
\begin{enumerate}
\item If $\gamma (s_0) \in T^{\ast }\L \cup \G_d$ then $\dot{\gamma}(s_{0})=H_{p_A}(\gamma(s_{0}))$.

\item If $\gamma (s_0) \in \G \setminus \G_d$ then $\dot{\gamma}(s_0)=H_{p_A}^G(\gamma(s_0))$, where  $H_{p_A}^G= H_{p_A} + (H_{p_A}^2x_{n} /H_{x_{n}}^2p_A)H_{x_{n}}$.

\item For every $s_0 \in B$, the two limits $\gamma(s_0 \pm 0)$ exist and are the two different points of the same hyperbolic fiber of the projection $\pi$.
\end{enumerate}
\end{defi}

Several remarks are in order:
\begin{enumerate}

\item If  $\Omega$ has no contact of infinite order  with its tangents, the Melrose-Sj\"ostrand flow is globally well defined, see \cite{MeSj}. 

\item In the interior, i.e., in $T^{*}\L$, a generalized bicharacteristic  is simply a classical bicharacteristic ray of the wave operator whose projection on the space variable is a geodesic of $\Omega$ equipped with the metric $(a^{ij})=(a_{ij})^{-1}$ .

\item  Finally, any generalized bicharacteristic ray $\gamma$ can be considered as a continuous map on the interval $I$ with values in $T^{\ast }_b\L$.
\end{enumerate}

\subsection{Sets of interest}\label{Subsec-Sets}

In what follows, we introduce several geometric sets associated with the Hamiltonian flow and linked to the observation region $\omega_{T}$. 
\medskip

\noindent {\bf A microlocal open subset of $T^{*}_{b} \mathcal{L}$.} We first introduce the set $\mathscr{R} (\omega_{T})$ defined by 
\begin{align}
	\mathscr{R} (\omega_{T}) 
	& 
	= 
	\{\rho = (t,x, \tau, \xi) \in T^{*}_{b} \mathcal{L} \backslash 0,  \text{ s. t. } \rho \notin  \Sigma_{A} \text{ or } \gamma_{\rho}(\R) \cap T^{*}_{b}\omega_{T} \neq \emptyset \}, 
\end{align}
which is the union of the set in which $P_{A}$ is elliptic, and of the set corresponding to the collection of bicharacteristic rays that meet the observation set $T^{*}_{b} \omega_{T}$. As we will see next, this is the set on which we can recover regularity properties on solutions of the wave equation from the regularity of the solution on $\omega_{T}$. 
\medskip

\noindent {\bf A microlocal open subset of $T^{*}_{b} \Omega$.} Another set, which will be of interest when discussing the recovery of microlocal information at $t = 0$, is the set $\mathscr{R}_{0}(\omega_{T})$ defined by
\begin{multline}\label{GCC-1}
\mathscr{R}_{0}(\omega_{T})
= \Big\{(x,\xi) \in T_{b}^{*}\Omega\backslash 0, 
\\
\text{ s. t. any $\gamma_\rho$}
\text{ emanating from $(x,\xi)$ at $t = 0$ satisfies } \gamma_{\rho} \cap T^{*}_{b}(\omega_{T}) \neq \emptyset \Big\}.
 \end{multline}
 In other words, $(x, \xi) \in \mathscr{R}_{0}(\omega_{T})$ if any bicharacteristic ray emanating from $(x, \xi)$ at $t=0$ enters in $\omega$ before the time $T$.
Let us emphasize immediately that for any $(x, \xi) \in \mathscr{R}_{0}(\omega_{T})$, there is at least two bicharacteristics emanating from $(x, \xi)$ at $t = 0$. 

To be more precise, we introduce the map $ j : T^{*}_{b}\L_{\vert t=0} \longrightarrow  T^{*}_{b}\Omega $  defined by
\begin{equation}
\left\{ 
\begin{array}{ll}
j(0,x; \tau,\xi) = (x,\xi) \quad \text{if} \quad (x,\xi) \in T^{*}\Omega ,
 \\
 j(0,x;\tau,\xi') = (x,\xi') \quad \text{if} \quad (x,\xi') \in T^{*}\d\Omega .
 \end{array}
  \right.  
	\label{map-j}
\end{equation}

In the sequel we will denote by $(x,\xi)$ the current point of $T^{*}_{b}\Omega $. If $x$ is a boundary point, $(x,\xi)$ has to be understood as $(x,\xi') \in T^{*}\d\Omega$, that is $\xi' \in \R^{n-1}$.

Recalling that  $\Sigma_{A} = \pi ( \text{Char}(P_{A}))$, we note that for $\tilde \rho = (x,\xi) \in T^{*}_{b}\Omega$, the set  $j^{-1}\{\tilde \rho\}\cap \Sigma_{A}$ is not empty.

Now, we make precise  the notion of  bicharacteristic curves of $P_{A}$, denoted by $\gamma$, emanating from  $(x, \xi)$ at $\{t=0\}$. Consider $\tilde \rho=(x,\xi) \in T_{b}^{*}\Omega\backslash 0$.
\begin{itemize}
\item If $x$ is an interior point, that is $x\in \Omega$, we have $j^{-1}\{\tilde \rho\}\cap \Sigma_{A} = \{(0,x;\tau= \pm\vert \xi\vert_{x}, \xi)\}$. Therefore, we have two  bicharacteristic curves issued from $\tilde \rho$, namely the  curve $\gamma^{+}$ issued from the point $\rho_{+}=(0,x;\tau=\vert \xi\vert_{x}, \xi)$,  and the curve $\gamma^{-}$ issued from the point $\rho_{-}=(0,x;\tau=-\vert \xi\vert_{x}, \xi)$.
\item If $x$ is a boundary  point, that is $x\in \d\Omega$, working in local geodesic coordinates, we have in this case $ \rho = (0,x;\tau,\xi') \in \Sigma_{A} \Leftrightarrow \tau^{2}\geq \vert \xi'\vert^{2}_{x}$.  

\textbf{a)} If $\vert\tau\vert  = \vert\xi'\vert_{x}$, we are dealing with a glancing point, and we know that for each $\tau= \pm\vert \xi'\vert_{x}$, there exists a unique ray  $\gamma_{\rho}$ passing through $\rho=(t=0,x; \tau, \xi')$. More precisely, if  $\rho \in \G_{d}$, we then identify $\rho$ to the point $(t=0,x; \tau= \pm\vert \xi'\vert_{x}, \xi', \xi_{n}=0) \in T^{*}\R^{d+1}$, and $\gamma_{\rho}$ is an integral curve of the (free) hamiltonian field $H_{p}$.  And if $\rho \in  \G\backslash\G_{d}$, then $\gamma_{\rho}$ is an integral curve  of the gliding field $H_{p}^{G}$, see Definition \ref{def:gene-bichar} .


\textbf{b)} If $\vert\tau\vert  > \vert\xi'\vert_{x}$, we are dealing with a hyperbolic point. $\gamma$ is then one of the two hyperbolic fibers of $P_{A}$ at $\rho$. According to the hamiltonian equations, on sees that    the  bicharacteristic curve  $\gamma$ corresponds to the integral curve of the (free) hamiltonian field $H_{p}$  issued from the point $ \rho_{-}=\big(0,x; - |\tau|,\xi', \xi_{n}=+\sqrt{\tau^{2} - \vert\xi'\vert_{x}^{2}} \, \big)$ or $ \rho_{+}=\big(0,x; + |\tau|,\xi', \xi_{n}=+\sqrt{\tau^{2} - \vert\xi'\vert_{x}^{2}} \, \big)$.
\end{itemize}

Let us emphasize that 
$$
	(x, \xi) \in \mathscr{R}_{0}(\omega_{T}) 
	\Leftrightarrow
	j^{-1}(x,\xi) \subset \mathscr{R}(\omega_{T}).
$$
\medskip

\noindent {\bf An open subset of $\overline{\Omega}$.} The last set of interest is the set 
$$
	\mathcal{O}(\omega_{T}) = 
	\{ x \in \overline\Omega, \text{ s. t. } T^{*}_{b}\Omega_{|\{x\}}\backslash 0 \subset \mathscr{R}_{0}(\omega_{T})
	\}, 
$$
which corresponds to the set of $x \in \overline\Omega$, from which all bicharacteristics emanating  at $t = 0$ meet the observation set $\omega_{T}$. 
\bigskip

Let us finally point out some basic properties of the sets $\mathscr{R}(\omega_{T})$, $\mathscr{R}_{0}(\omega_{T})$, and $\mathcal{O}(\omega_{T})$ :
\begin{itemize}
\item All these sets are non-empty since obviously $T^*_b\L_{\vert \omega_{T}}\backslash 0 \subset \mathscr{R}(\omega_{T})$, $T^{*}_b\Omega_{\vert \omega}\backslash 0 \subset \mathscr{R}_0(\omega_{T})$, and $\omega \subset \mathcal{O}(\omega_{T})$.  
\item   $\mathscr{R}(\omega_{T})$, $\mathscr{R}_{0}(\omega_{T})$, and $\mathcal{O}(\omega_{T})$ respectively are open subsets of  $T_{b}^{*}\mathcal{L}$, $T_{b}^{*}\Omega$ and $\overline{\Omega}$, according to the continuity of the Melrose-Sj\"ostrand flow.  
\end{itemize}

Note that the classical GCC for $\omega_{T}$ can be simply expressed as one of the following equivalent formulations: $\mathscr{R}(\omega_{T}) = T_{b}^* \mathcal{L} \backslash 0$, 
$\mathcal{R}_{0}(\omega_{T})  = T_{b}^{*}\Omega \backslash 0$, or $\mathcal{O}(\omega_{T}) = \overline\Omega$.

The geometric condition of Theorem \ref{Theo-Obs-1} can in fact be simply stated as $\overline{\mathscr{O}} \subset \mathcal{O}(\omega_{T})$. In other words, Theorem \ref{Theo-Obs-1} applies  for any open set $\mathscr{O}$ strictly included in $\mathcal{O}(\omega_{T})$.

\subsection{Pseudo-differential operators}\label{pdo's}

Following \cite{Lebeau96} ,we define the set $\mathcal{A}$ of pseudo-differential operators on $\R \times \R^{n}$ of the form $Q=Q_{i}+Q_{\d}$ where $Q_{i}$ is a classical pseudo-differential operator,  compactly supported in $\L$ and $Q_{\d}$ is a classical   pseudo-differential operator tangential to the boundary $\d\L$, compactly supported near $\d\L$. More precisely, $Q_{i}=\varphi Q_{i}\varphi$ for some $\varphi \in \mathscr{C}^{\infty}_{0}(\L)$,  and $Q_{\d}=\psi Q_{\d}\psi$ for some $\psi \in \mathscr{C}_{0}^{\infty}(U_{\d\L})$, where $U_{\d\L}$  is a small neighborhood of $\d\L$ in $\R \times \R^{n}$. For $s \in \R$, $\mathcal{A}^{s}$ denotes  the set of elements of order $s$ of $\mathcal{A}$.

 In a similar way, we also define the set 
$\mathcal{B}$ of pseudo-differential operators on $ \R^{n}$, i.e., on the space variable,  of the form $\psi=\psi_{i}+\psi_{\d}$ where $\psi_{i}$ is a classical pseudo-differential operator,  compactly supported in $\Omega$, and $\psi_{\d}$ is a classical  pseudo-differential operator tangential to the boundary $\d\Omega$, compactly supported near $\d\Omega$. Similarly as above, for $s \in \R$, $\mathcal{B}^{s}$ denotes the set of the elements of $\mathcal{B}$ of order $s$.

\subsection{Wave front sets and propagation results}\label{wave front}
In this section, we recall the notion of wave front set  up to the boundary and the classical propagation results of Melrose--Sj\"ostrand \cite{MeSj} and H\"ormander \cite{Hormander-III}. 

For a distribution $u$ defined on the cylinder $ \mathcal{L} = \mathbb{R} \times \Omega$, we define the $H^s$-wave front set up to the boundary, denoted $WF^s_b(u)$, as a subset of the compressed cotangent bundle in the sense of Melrose–Sjöstrand,
$T^*_b\L=T^{\ast }\L\cup T^{\ast }\partial \L$.
This set coincides with the classical wave front set $WF^s(u)$ in the interior of $\mathcal{L}$, i.e., in $T^* \mathcal{L}$, and extends the notion to describe the $H^s$-microlocal regularity of $u$ up to the boundary.
We follow here the definition of Chazarain (see \cite{Chazarain}) that is used by Melrose and Sj\"ostrand (see  \cite{MeSj}). In addition, for solutions of $P_{A}u\in \mathscr{C}^{\infty}$, it agrees with the intrinsic notion of Melrose  (see H\"ormander \cite[Cor.18.3.33]{Hormander-III}), which does not depend on $P_{A}$.
In the following, we use the spaces of pseudodifferential operators  $\mathcal{A}^{0}$ and $\mathcal{B}^{0}$, introduced in Section \ref{pdo's}.

Consider $s \in \R$ and $u \in \mathscr{D}'(\L)$ (later, $u$ will be a solution of $P_{A}u=0$ in $\L$). Also, for $q \in \R^{n+1}$ and $r>0$,  denote by $B_{r}(q)$ the Euclidean ball of center $q$ and radius $r$.

\begin{defi}\label{wave front-1}
For $\rho=(q,\eta) \in T_{b}^{\ast }\L$, we say that $\rho \notin WF^{s}_{b}u$ if there exists a  pseudodifferential operator $Q \in \mathcal{A}^{0}$  such that $Q$ is elliptic at $\rho$, and $Qu \in H^{s}(B_{r}(q)\cap\L)$ for some $r>0$.

More precisely, 
\begin{itemize}
\item  If $\rho=(q,\eta) \in T^{\ast }\L$, i.e., $q$ is an interior point,  there exists a  pseudodifferential operator $Q=Q_{i} \in \mathcal{A}^{0}$, elliptic at $\rho$, such that $Qu \in H^{s}(B_{r}(q))$ for some $r>0$ ,  $B_{r}(q) \subset \L$ .
\item If $\rho=(q,\eta) \in T^{\ast }\d\L$, i.e., $q$ is a boundary point,  there exists a  tangential pseudodifferential operator $Q=Q_{\d} \in \mathcal{A}^{0}$, elliptic at $\rho$, such that $Qu \in H^{s}(B_{r}(q)\cap \L)$ for some $r>0$.
\end{itemize}

\end{defi}
\begin{rema}
For $v \in \mathcal{D}'(\Omega)$ and $\rho' \in T_{b}^{\ast }\Omega$, we have a similar definition for  $\rho' \notin WF^{s}_{b}v$.
\end{rema}

Here we recall the Melrose-Sj\"ostrand theorem for propagation of regularity. For the convenience of the reader, we give a statement adapted to the framework of system \eqref{waveequation}.

Remind that for $\rho \in \Sigma_{A} \subset T_{b}^{\ast }\L$, we denote by $\gamma_{\rho}$ the generalized bicharacteristic curve of $P_{A}$, issued from $\rho$ as described in Section \ref{Subsec-Sets} above.

\begin{theorem}[Melrose-Sj\"ostrand \cite{MeSj}]\label{M-S. theorem}
Let $u$ be a solution of system \eqref{waveequation} with $(u_{0},u_{1})\in L^{2}(\Omega)\times H^{-1}(\Omega)$, and assume that a point $\rho \in T_{b}^{\ast }\L$ is such that $\rho \notin WF^{1}_{b}u$. Then $\gamma_{\rho}\cap WF^{1}_{b}u = \emptyset$. 
\end{theorem}

\subsection{Proof of Lemma \ref{regularity}}\label{Subsec-Proof-Lemma-reg}

First, we notice that  if $(u_{0},u_{1})\in L^2(\Omega)\times H^{-1}(\Omega)$, and the corresponding solution satisfies $\d_{t}u \in L^{2}(\omega_{T})$, then $u$ lies in $H^{1}_{loc}(\omega_{T})$ by microlocal elliptic regularity.

We will deduce that  $u$ actually belongs to $ H^{1}((0,T)\times \Omega)$.
Indeed, let us consider $\rho_{0} \in T^{*}_{b}((0,T)\times \Omega)$. If  $\rho_{0}$ is an elliptic point (independently if it is an interior or a boundary point), it is classical that $\rho_{0}\notin WF_{b}^{2}u$, i.e., $u$ is in $H^{2}$  microlocally near $\rho_{0}$. Here, a special care must be taken at the boundary,  see H\"ormander \cite [Theorem 20.1.14]{Hormander-III} .

If  $\rho_{0}$ is not an elliptic point, denote by $\gamma_{\rho_{0}}$ the generalized bicharacteristic ray issued from $\rho_{0}$. We have two possibilities : a) $\gamma_{\rho_{0}}$ intersects $T^{*}(\omega_{T})$, or  b) $\gamma_{\rho_{0}} \cap T^{*}(\omega_{T}) = \emptyset$.

In case a), since $u \in H^{1}_{loc}(\omega_{T})$, $\rho_{0} \notin WF_{b}^{1}u$ by propagation of the $H^{1}$-wave front, see Theorem \ref{M-S. theorem} above.

In case b), following $\gamma_{\rho_{0}}$ backward in time,  let us set $\rho_{1}=\gamma_{\rho_{0}} \cap \{t=0\}$. According to the microlocal assumption on the pair $(\omega,  \mathscr{O})$, we have $x(\rho_{1}) \notin \overline{\mathscr{O}}$. Therefore, the initial data $(u_{0},u_{1})$ is vanishing in a neighborhood of $x(\rho_{1})$ and so does the solution $u$ in some space-time cylinder $(-\alpha, \alpha)\times B(x(\rho_{1}),r)$, $\alpha >0, r>0$ small. Consequently, $\rho_{1} \notin WF_{b}^{1}u$, and again, we obtain $\rho_{0} \notin WF_{b}^{1}u$, by propagation of the $H^{1}$-wave front up to the boundary.

Accordingly, the solution $u$ of the wave system \eqref{waveequation} lies in $H^{1}((0,T)\times \Omega)$. Thus, by conservation of energy in time, we also conclude that the initial data has finite energy, i.e.,$(u_0, u_1)$ belongs to $H^1_0 (\Omega)\times L^2(\Omega)$. \hfill $\square$

\begin{rema}
	\label{Rema-Lem-regularity-boundary}
	The proof of the propagation of regularity from the boundary stated in Lemma \ref{regularity-boundary} follows the same strategy as the one of Lemma \ref{regularity}. The only difference is that one has to use the  propagation of the $H^1$-wave front set from an observation on the boundary, namely out of the information that $\d_{n}u_{\vert \d\Omega} \in L^{2}(\Gamma_{T})$. This can be done on nondiffractive points.  Indeed, under the assumption $\d_{n}u_{\vert \d\Omega} \in L^{2}(\Gamma_{T})$, by the lifting lemma in \cite[Theorem 2.2]{B-L-R}, we deduce that every nondiffractive  point $\rho_{0} \in T^{*}(\Gamma_{T})$ satisfies $\rho_{0} \notin WF_{b}^{1}(u)$. In other words, the solution $u$ is $H^{1}$ microlocally near this point, up to the boundary. 
	This suffices to conclude the proof of Lemma \ref{regularity-boundary} in view of the imposed microlocal geometric condition on the pair $(\Gamma, \mathscr{O})$.%
\end{rema}

\subsection{Further technical results}

The goal of this section is to prove the following generalization of Lemma \ref{regularity}, which underlines the role played by the various sets $\mathscr{R}(\omega_{T})$, $\mathscr{R}_0(\omega_{T})$, and $\mathcal{O}(\omega_{T})$ introduced in Section \ref{Subsec-Sets}. 

\begin{lemma}\label{regularity-general}
	Let $u$ be a solution $u$ of \eqref{waveequation} with initial data $(u_{0},u_{1}) \in L^2(\Omega) \times H^{-1}(\Omega)$,  and  satisfying $\partial_t u \in L^{2}(\omega_{T} )$. 
Then  
	\begin{enumerate}
		\item[(1a)] $WF_b^1(u) \cap \mathscr{R}(\omega_{T}) = \emptyset$,
		\item[(2a)] $(WF_b^1(u_0) \cup WF_b^{0}(u_1)) \cap\mathscr{R}_0(\omega_{T}) = \emptyset$, 
		\item[(3a)] $(u_0, u_1)   \in H^1_{loc}(\mathcal{O}(\omega_{T})) \times L^{2}_{loc}(\mathcal{O}(\omega_{T}))$.
	\end{enumerate}

	Similarly, if $u$  solution of \eqref{waveequation}, with initial data $(u_{0},u_{1}) \in H^{-1}(\Omega)\times (H^{2}\cap H^1_0)'(\Omega)$, satisfies $u \in L^{2}(\omega_{T} )$, we have 
	\begin{enumerate}
		\item[(1b)] $WF_b^0(u) \cap \mathscr{R}(\omega_{T}) = \emptyset$,
		\item[(2b)] $(WF_b^0(u_0) \cup WF_b^{-1}(u_1)) \cap\mathscr{R}_0(\omega_{T}) = \emptyset$, 
		\item[(3b)] $(u_0, u_1)   \in L^2_{loc}(\mathcal{O}(\omega_{T})) \times H^{-1}_{loc}(\mathcal{O}(\omega_{T}))$. 
	\end{enumerate}
\end{lemma}

\begin{proof}[Proof of Lemma \ref{regularity-general}]

\noindent

\noindent\emph{Item (1):}	The proof of this item is in fact included in the proof of Lemma \ref{regularity}. 
	Let us consider $u$ a solution of \eqref{waveequation} with initial data $(u_{0},u_{1}) \in L^2(\Omega) \times H^{-1}(\Omega)$,  and  satisfying $\partial_t u \in L^{2}(\omega_{T} )$ (the other case, being completely similar, it is left to the reader).
	 
	For $\rho_0 \in T_b^*\mathcal{L} \cap  \mathscr{R}(\omega_{T})$, there are two possibilities: either $\rho_0 \notin \Sigma_A$, or $\rho_0 \in \Sigma_A$ and $\gamma_{\rho_0} \cap T^*_b(\omega_{T}) \neq \emptyset$: When $\rho_0 \notin \Sigma_A$, the operator is elliptic, so $u$ is in $ H^2$ microlocally near $\rho_0$; when $\rho_0 \in \Sigma_A$ and $\gamma_{\rho_0} \cap T^*_b(\omega_{T}) \neq \emptyset$, Theorem \ref{M-S. theorem} guarantees that $u$ is in $H^1$ microlocally near $\rho_0$.

\noindent\emph{Item (3):}
	Here again, we consider a solution $u$  of \eqref{waveequation} with initial data $(u_{0},u_{1}) \in L^2(\Omega) \times H^{-1}(\Omega)$,  and  satisfying $\partial_t u \in L^{2}(\omega_{T} )$ (the other case being completely similar, it is left to the reader).
	
	Let $x_0 \in \mathcal{O}(\omega_{T})$. Since $\mathcal{O}(\omega_{T})$ is open in $\overline\Omega$, by continuity of the flow of Melrose-Sj\"ostrand, there exists $\varepsilon>0$ and $\varphi \in \mathscr{C}^\infty_{c}((-\varepsilon, \varepsilon) \times \overline\Omega)$ such that $\varphi \equiv 1$ in a neighborhood of $(0, x_0)$, such that all rays emanating from the support of $\varphi$ meet $\omega_{T}$.
	
	By item (1), we know that $u$ is locally $H^1$ near any point of the support of $\varphi$. To conclude, we simply note that $v = \varphi u$ solves 
	\begin{equation*}
\left\{ 
\begin{array}{ll}
	P_{A}v =  [P_A,\varphi]u \in L^2(\mathcal{L}) 
		\\
v_{\vert\d\Omega} = 0
	\\
	 v(\varepsilon, \cdot)= \d_{t}v(\varepsilon, \cdot)=0
  \end{array}
  \right.  
	\end{equation*}
	and thus $(v, \partial_t v) |_{t = 0} \in H^1_0(\Omega) \times L^2(\Omega)$. This concludes the proof since $u = v$ in a neighborhood of $x_0$.

\noindent\emph{Item (2):}
The proof of item (2) of Lemma \ref{regularity-general} is more involved and requires several additional technical results, which are presented in detail below.

Using the continuity of Melrose-Sj\"ostrand flow, there exists $\varepsilon >0$ such that the set 
\begin{align*} 
E_{\varepsilon}=\{ \rho = (t,x;\tau,\xi) \in T^{*}_{b}\L, \,\,(x,\xi) \in \mathscr{R}_{0}(\omega_{T}) , \, \,  t \in (0, \varepsilon ) \} 
\end{align*}   
satisfies $E_{\varepsilon} \subset \mathscr{R}(\omega_{T})$ and thus $E_\varepsilon \cap WF^{1}_{b}u = \emptyset $ by item (1) of Lemma \ref{regularity-general}. 

Our first goal is to prove that for any $\psi (x,D_{x}) \in \mathcal{B}^{0}$ supported in $\mathscr{R}_{0}(\omega_{T})$, we have $\psi (x, D_x) u \in H^1((0,\varepsilon) \times \Omega)$. And for this end, we will prove that any point $\rho_{0} =(t_0,x_{0};\tau_{0},\xi_{0}) \in T^{*}_{b}\L, t_{0} \in (0 ,\varepsilon )$, satisfies $\rho_{0} \notin WF^{1}_{b}(\psi (x, D_x)u)$.

\smallskip

\noindent {\it Case 1:} $\rho_{0} \in T^{*}\L$, i.e., it is an interior point and $x_{0} \in \Omega$.  Take $\varphi = \varphi (t,x) \in \mathscr{C}_{0}^{\infty}((0 ,\varepsilon )\times \Omega)$, supported near $(T,x_{0})$.
In the operators space $\mathcal{A}^{0}$ consider the local identity partition 
\begin{align*} 
\varphi (t,x) = q_{1}(t,x;D_{t},D_{x})+ q_{2}(t,x;D_{t},D_{x}) +\mathcal{R}, 
\end{align*} 
with 
\begin{align*} 
\supp (q_{1}) \subset \left\{\frac{1}{2} \vert\xi \vert \leq \vert\tau \vert \leq  2 \vert\xi \vert \right\} ,
\\
\supp (q_{2}) \subset \left\{\vert\tau \vert \leq \frac{3}{4} \vert\xi \vert   \right\} \cup \left\{ \vert\tau \vert \geq   \frac{3}{2} \vert\xi \vert \right\},
\end{align*} 
and $\supp_{(t,x)}(q_{j})$, $j=1,2$ is a compact of $(0 ,\varepsilon )\times \Omega$, and $\mathcal{R}$ is infinitely smoothing.

Clearly, $\supp (q_{2})$ is contained in the elliptic set of $T^{*}\L$, so $q_{2}(t,x;D_{t},D_{x})u \in H^{1}(\L_{\varepsilon})$, and  $\psi(x,D_{x})q_{2}u  \in H^{1}(\L_{\varepsilon})$ .

Let us now examine the first term $q_{1}(t,x;D_{t},D_{x}) u$. Here we notice that the composition $\psi(x,D_{x}) q_{1}(t,x;D_{t},D_{x})$  provides a well defined global pseudodifferential operator, according to \cite[Th. 18.1.35]{Hormander-III}. In addition, if $\rho = (t,x;\tau,\xi) \in \supp(\sigma(\psi q_{1}))$, where $\sigma$ denotes the symbol of the operator $\psi q_1$, $\rho$ is either an elliptic or hyperbolic point lying in the set $E_{\varepsilon}$, which does not intersect $WF^1u$. Therefore, $\psi(x,D_{x}) q_{1}(t,x;D_{t},D_{x}) u \in H^{1}(\L_{\varepsilon})$.

Hence  we deduce  $ \varphi \psi(x,D_{x})u= \psi(x,D_{x})(\varphi u) -[\varphi, \psi (x,D_{x})]u  \in H^{1}(\L_{\varepsilon})$.

\smallskip
\noindent {\it Case 2:} $\rho_{0} \in T^{*}\d\L$, i.e., it is a boundary point and $x_{0} \in \d\Omega$.

Here we shall work in a system of local geodesic coordinates $(t,x',x_{n};\tau,\xi',\xi_{n})$ with $\d\Omega = \{x_{n}=0 \}$ and $\Omega = \{x_{n}>0 \}$, see Section \ref{geo}. Hence we will set $\rho_{0}=(T,x'_{0},\tau_{0},\xi'_{0})$. Recall that the operators  of $\mathcal{A}^{0}$ (resp. of $\mathcal{B}{0}$) take the form $q(x_{n},t,x',D_{t},D_{x'})$  (resp. $\psi(x_{n},x',D_{x'})$).

As in Case 1 above, we consider $\varphi = \varphi (t,x) \in \mathscr{C}_{0}^{\infty}((0 ,\varepsilon )\times \R^{n})$, supported near $(T,x_{0})$, and a local partition of the identity with tangential pseudodifferential operators,  of the form 
\begin{align*} 
\varphi (t,x) = q_{1}(x_{n}, t,x';D_{t},D_{x'})+ q_{2}(x_{n},t,x';D_{t},D_{x'}) + \mathcal{R}, 
\end{align*} 
with, this time, 
\begin{align*} 
\supp (q_{1}) \subset \{ \vert\tau \vert \leq 3 \vert\xi' \vert  \} \quad  \text{and } \quad \supp (q_{2}) \subset \{\vert\tau \vert  \geq 2\vert\xi' \vert \} .
\end{align*} 
With notations of Section \ref{boundary-geometry}, $\supp (q_{2}) \subset \H$, the hyperbolic subset of $T^{*}\d\L$. Therefore since $\supp(\psi) \subset \mathscr{R}_{0}(\omega_{T})$, we get $ q_{2}(x_{n},t,x';D_{t},D_{x'})\psi(x_{n},x', D_{x'}) u \in H^{1}(\L_{T})$ .

Finally, as in Case 1, we notice that the composition $q_{1}(x_{n},t,x';D_{t},D_{x'})\psi(x_{n},x', D_{x'})$  provides a well defined global   tangential pseudodifferential operator, see \cite[Th. 18.1.35]{Hormander-III}, whose support is contained in $\E \cup \mathscr{R}(\omega_{T})$ . And this yields $q_{1}(x_{n},t,x';D_{t},D_{x'})\psi(x_{n},x', D_{x'}) u \in H^{1}(\L_{T})$ according to item (1).
\medskip

Let us now examine the regularity of the trace $(u_{0},u_{1})$. For this, consider  a function $h(t) \in \mathscr{C}^{\infty}_{0}(\R)$, $h(t)=1$ for $\vert t \vert \leq \varepsilon/2$ and $h(t)=0$ for $\vert t \vert \geq 3/4\varepsilon$. For $u(t,x)$ solution to \eqref{waveequation}, the function $v = h(t)\psi(x,D_{x})u$ satisfies the wave system
\begin{equation}
\left\{ 
\begin{array}{ll}
P_{A}v =  [\d_{t}^{2},h(t)]\psi(x,D_{x})u - h(t)[\Delta_{A},\psi(x,D_{x})]u , 
\\
v_{\vert\d\Omega} = 0
\\
 v(\varepsilon, \cdot)= \d_{t}v(\varepsilon, \cdot)=0.
  \end{array}
  \right.  
\end{equation}

The right hand side of this equation lies in $L^{2}((0,\varepsilon)\times \Omega)$ according to the argument above. Therefore $(v(0, \cdot), \partial_t v(0, \cdot))=\psi(x,D_{x})(u_{0},u_{1})$ belongs to $H^1_0 (\Omega)\times L^2(\Omega)$. This ends the proof of  item (2) in Lemma \ref{regularity-general}, since $\psi(x, D)$ is any operator in $\mathcal{B}^{0}$ supported in $\mathscr{R}_{0}(\omega_{T})$.
\end{proof}

\subsection{A $1$-d example} 

To better understand the geometric statements given by Lemma \ref{regularity-general}, we briefly study the $1$-d case when $\Omega =( -10,10)$.

In this case, the wave equation 
\begin{equation}
\left\{ 
\begin{array}{ll}
	\partial _{t}^{2}u-\partial_{x}^{2}u=0,\quad &\text{in } \,\,(0,T) \times (-10, 10), 
	\\ 
	u(t,-10) = u(t, 10) = 0,\quad &\text{on } \,\,(0,T),
	\\
	(u(0,\cdot),\partial _{t}u(0,\cdot ))=(u_{0},u_{1}), & 
	\end{array}%
	\right.  
	\label{waveequation$1$-d}
\end{equation}
can be solved explicitly using the characteristics. Indeed, setting 
$$
	w_{+}(t,x) = (\partial_{t} u - \partial_{x} u)(t,x), 
	\qquad
	w_{-}(t,x) = (\partial_{t} u + \partial_{x} u)(t,x), 
	\qquad 
	\text{ for } (t,x) \in (0,T) \times (-10,10), 
$$
the $1$-d wave equation can be recast into a system of transport equation coupled from the boundary
\begin{equation}
\left\{ 
\begin{array}{ll}
	\partial _{t} w_{+}+ \partial_{x}w_{+} =0,\quad &\text{in } \,\,(0,T) \times (-10, 10), 
	\\
	\partial _{t} w_{-}- \partial_{x}w_{-} =0,\quad &\text{in } \,\,(0,T) \times (-10, 10), 
	\\ 
	(w_{+}+w_{-})(t,-10) = (w_{+}+w_{-}) (t, 10) = 0,\quad &\text{on } \,\,(0,T),
	\\
	(w_{+}(0,\cdot),w_{-}(0,\cdot ))=(u_{1} - \partial_{x }u_{0},u_{1} + \partial_{x }u_{0}), & 
	\end{array}%
	\right.  
	\label{waveequation$1$-d-two-transport}
\end{equation}
In this case, the bicharacteristic rays are particularly simple: they are the curves $t \mapsto x_{0} +t$ and $t \mapsto x_{0} -t$ for $x_{0} \in (-10,10)$ while these curves stay in the domain, bouncing back when meeting the boundary. Accordingly, in $1$-d, we can identify the characteristic manifold $\text{Char}(P_{A})$ with $\R \times \Omega \times \{-1, 1\}$, depending if $\tau = |\xi|$, corresponding to $\epsilon = 1$, or $\tau = - |\xi|$ corresponding to $\epsilon = -1$.

Let us now fix  $\omega = (-2, -1) \cup (1, 2)$, and $T = 3$. The sets $\mathscr{R}(\omega_{T})$, $\mathscr{R}_{0}(\omega_{T})$ and $\mathcal{O}(\omega_{T})$ can then computed explicitly:
\begin{equation}
	\mathscr{R}(\omega_{T})|_{t \in [0,T]} 
	= 
	\mathscr{R}(\omega_{T})^{+} \cup \mathscr{R}(\omega_{T})^{-}
\end{equation}
with 
\begin{align*}
	\mathscr{R}(\omega_{T})^{+} &= 
	\{ (t,x, +1) \text{ s.t. } t \in [0,T] \text{ and } -5+t <x < 2+t \}
	\\
	\mathscr{R}(\omega_{T})^{-} &= 
	\{ (t,x, -1) \text{ s.t. } t \in [0,T] \text{ and } -2-t <x < 5-t \},
\end{align*}
and 
\begin{align*}
	& \mathscr{R}_{0}(\omega_{T})  = (-2, 2) \times \R^{*},
	\\
	& \mathcal{O}(\omega_{T}) = (-2, 2).
\end{align*}
To illustrate Lemma \ref{regularity-general}, due to the structure of the solutions of the wave given by \eqref{waveequation$1$-d-two-transport}, it is clear that if $(u_{0}, u_{1}) \in L^2(\Omega) \times H^{-1}(\Omega)$ with $ u \in H^1(\omega_{T})$, which of course implies $w_{+}$ and $w_{-}$ belong to $L^2(\omega_{T})$, 
\begin{enumerate}
	\item $w_{+} \in L^{2}(\mathscr{R}(\omega_{T})^{+})$ and $w_{-}  \in L^{2}(\mathscr{R}(\omega_{T})^{+})$, 
	\item $w_{+}|_{t = 0} \in L^2(-5,2)$ and $w_{-}|_{t = 0} \in L^2(-2,5)$, 
	\item $(u_{0}, u_{1}) \in H^1(-2,2) \times L^{2}(-2,2)$.
\end{enumerate}
It is also clear due to the explicit character of the solutions of \eqref{waveequation$1$-d-two-transport} that we cannot improve these sets of regularity for general data.

\section{Further observability results}\label{sec-further results}

The aim of this section is to refine Theorem \ref{Theo-Obs-1} by analysing which microlocal components of the initial data of general wave solutions can be effectively observed from measurements taken on $\omega_T$.

Up to this point, our focus has been on initial data supported in a set \( \overline{\mathscr{O}} \), such that the pair \( (\omega, \mathscr{O}) \) satisfies the required microlocal geometric condition. We now adopt a complementary viewpoint: given a fixed observation region $\omega_T$, we seek to extract the maximum amount of information possible from the available measurements. As we shall see, we can recover, in a precise sense, the energy associated with the microlocal projection of the initial data that propagates along rays entering the observation region $\omega$.

The proofs of these refined results follow the same general strategy as before, relying in particular on Lemma \ref{regularity-general}, which ensures propagation of microlocal regularity. The final observability estimates depend on whether a suitable unique continuation property is available, which determines our ability to eliminate the compact remainder term.

\subsection{Statement of the results}\label{statement-results}

We start with the following microlocal observability estimates.

\begin{theorem}[Relaxed microlocal observability estimates] \label{R-obs}
	Let $\omega$ be a non-empty open set of $\Omega$ and $T >0$. 
		\begin{enumerate}
		\item[(1a)] For every operator $\psi (t,x,D_{t},D_{x}) \in \mathcal{A}^{0}$ such that $\supp(\psi)\cap T^{*}_{b}\L \subset \mathcal{R}(\omega_{T})$, there exists $C>0$ such that for every initial data $(u_{0},u_{1}) \in H^1_0(\Omega) \times L^{2}(\Omega)$, 
the solution of  \eqref{waveequation} satisfies 
\begin{equation}\label{R-obs-t-x-H1-micro}
	\Vert \psi (t,x,D_{t},D_{x})u \Vert_{H^{1}(\L_{T})} \leq C \Vert \partial_t u\Vert_{L^{2}(\omega_{T} )} +C \Vert u \Vert_{L^{2}(\L_{T})}.
\end{equation}
		\item[(2a)] For every operator $\psi (x,D_{x}) \in \mathcal{B}^{0}$ such that $\supp(\psi)\cap T^{*}_{b}\Omega \subset \mathcal{R}_{0}(\omega_{T})$, there exists $C>0$ such that  for every initial data $(u_{0},u_{1}) \in H^1_0(\Omega) \times L^{2}(\Omega)$, 
the solution of  \eqref{waveequation} satisfies 
\begin{equation}\label{R-obs-H1-micro}
\Vert \psi (x,D_{x})(u_{0},u_{1})\Vert_{H^1_0(\Omega) \times L^{2}(\Omega)} \leq C \Vert \partial_t u\Vert_{L^{2}(\omega_{T} )} +C \Vert(u_{0},u_{1})\Vert_{L^2(\Omega) \times H^{-1}(\Omega)}.
\end{equation}
		\item[(3a)] For every function $\psi= \psi (x) \in \mathscr{C}^\infty_{c}(\mathcal{O}(\omega_{T}))$, there exists $C>0$ such that for every initial data $(u_{0},u_{1}) \in H^1_0(\Omega) \times L^{2}(\Omega)$, 
the solution of  \eqref{waveequation} satisfies 
\begin{equation}\label{R-obs-H1-local}
\Vert \psi (x)(u_{0},u_{1})\Vert_{H^1_0(\Omega) \times L^{2}(\Omega)} \leq C \Vert \partial_t u\Vert_{L^{2}(\omega_{T} )} +C \Vert(u_{0},u_{1})\Vert_{L^2(\Omega) \times H^{-1}(\Omega)}.
\end{equation}
\end{enumerate}

Similarly,
	
	\begin{enumerate}
		\item[(1b)] For every operator $\psi (t,x,D_{t},D_{x}) \in \mathcal{A}^{0}$ such that $\supp(\psi)\cap T^{*}_{b}\L \subset \mathcal{R}(\omega_{T})$, there exists $C>0$ such that for every initial data $(u_{0},u_{1}) \in L^{2}(\Omega)\times H^{-1}(\Omega)$, 
the solution of  \eqref{waveequation} satisfies 
\begin{equation}\label{R-obs-t-x-L2-micro}
	\Vert \psi (t,x,D_{t},D_{x})u \Vert_{L^{2}(\L_{T})} \leq C \Vert u\Vert_{L^{2}(\omega_{T} )} +C \Vert u \Vert_{H^{-1}(\L_{T})}.
\end{equation}
		\item[(2b)] For every operator $\psi (x,D_{x}) \in \mathcal{B}^{0}$ such that $\supp(\psi)\cap T^{*}_{b}\Omega \subset \mathcal{R}_{0}(\omega_{T})$, there exists $C>0$ such that  for every initial data $(u_{0},u_{1}) \in L^{2}(\Omega)\times H^{-1}(\Omega)$, 
the solution of  \eqref{waveequation} satisfies 
\begin{equation}\label{R-obs-L2-micro}
\Vert \psi (x,D_{x})(u_{0},u_{1})\Vert_{L^{2}(\Omega)\times H^{-1}(\Omega)} \leq C \Vert u\Vert_{L^{2}(\omega_{T} )} +C \Vert(u_{0},u_{1})\Vert_{H^{-1}(\Omega)\times (H^{2}\cap H^1_0))'(\Omega)}.
\end{equation}
		\item[(3b)] For every function $\psi = \psi (x) \in \mathscr{C}^\infty_{c}(\mathcal{O}(\omega_{T}))$, there exists $C>0$ such that for every initial data $(u_{0},u_{1}) \in L^{2}(\Omega)\times H^{-1}(\Omega)$, 
the solution of  \eqref{waveequation} satisfies 
\begin{equation}\label{R-obs-L2-local}
\Vert \psi (x)(u_{0},u_{1})\Vert_{L^{2}(\Omega)\times H^{-1}(\Omega)} \leq C \Vert u\Vert_{L^{2}(\omega_{T} )} +C \Vert(u_{0},u_{1})\Vert_{H^{-1}(\Omega)\times (H^{2}\cap H^1_0))'(\Omega)}.
\end{equation}
\end{enumerate}
\end{theorem} 

\begin{rema}\label{Rk-Better-Local-Obs}
	Theorem \ref{R-obs}, items (3a) and (3b), and Lemma \ref{UC} allows to generalize the result of Theorem \ref{Theo-Obs-1} as follows:
	\begin{corollary} \label{Cor-Obs-1-bis}
		 For every function $\psi = \psi (x) \in \mathscr{C}^\infty_{c}(\mathcal{O}(\omega_{T}))$, there exists $C >0$ such that for any $(u_{0},u_{1})\in  H^1_0(\Omega) \times L^{2}(\Omega)$, the solution $u$ of \eqref{waveequation} satisfies the observability estimate 
	\begin{equation}
		\label{Distributed-Obs-Best-Geo-O-omega-H1}
		\Vert \psi (x)(u_{0},u_{1})\Vert_{ H^1_0(\Omega) \times L^{2}(\Omega)} 
		\leq C\Vert \partial_t u\Vert_{L^{2}(\omega_{T} )} +C \Vert  (1- \psi (x)) (u_{0},u_{1})\Vert_{L^{2}(\Omega) \times H^{-1}(\Omega)}.
	\end{equation}
	Similarly, for  any initial data $(u_{0},u_{1}) \in L^{2}(\Omega)\times H^{-1}(\Omega)$,
	 	\begin{equation}
		\label{Distributed-Obs-Best-Geo-O-omega-L2}
		\Vert \psi (x)(u_{0},u_{1})\Vert_{ L^{2}(\Omega) \times H^{-1}(\Omega)} 
		\leq C\Vert u\Vert_{L^{2}(\omega_{T} )} +C \Vert  (1- \psi (x)) (u_{0},u_{1})\Vert_{ H^{-1}(\Omega) \times (H^2 \cap H^1_0)'(\Omega)}. 
	\end{equation}
	\end{corollary} 

	The proof can be done similarly as the one of Theorem \ref{Theo-Obs-1}, passing from the estimates \eqref{R-obs-H1-local} and \eqref{R-obs-L2-local} to the estimates \eqref{Distributed-Obs-Best-Geo-O-omega-H1}--\eqref{Distributed-Obs-Best-Geo-O-omega-L2} by contradiction, and using the unique continuation result given by Lemma \ref{UC} for solutions of the wave equation with supported in $\supp \varphi$. Details are left to the reader. 
\end{rema}

\begin{rema}
Several comments are in order:
\begin{itemize}
\item The estimates in Theorem \ref{R-obs} hold true independently of any unique continuation consideration. They only rely on the propagation of regularity for solutions to the wave operator from $\omega_{T}$ to capture the energy of the projections of the  data determined in the sets $\mathscr{R}(\omega_{T})$, $\mathscr{R}_0(\omega_{T})$ and $\mathcal{O}(\omega_{T})$ introduced in Section \ref{Subsec-Sets} according to Lemma \ref{regularity-general}. 

\item 
The constants appearing in Theorem \ref{R-obs} depend on both the observation time $T$ and the underlying metric $(a^{ij}(x))$. The dependence on the metric is 
very subtle  and difficult to quantify explicitly, as it involves a large number of derivatives of the coefficients $a^{ij}$. This dependence can, in principle, be traced through a careful analysis of the proof of Theorem \ref{M-S. theorem} in \cite{MeSj}. However, it is not readily expressible in closed form, since the operator norm of a pseudo-differential operator typically involves multiple derivatives of its symbol, making the exact dependence implicit and technically intricate.

\item Due to the potential failure of the unique continuation property, an additive remainder term is required to ensure the validity of the inequalities. In the following, we will discuss how this remainder can be weakened or even eliminated under additional geometric assumptions.
\end{itemize}
\end{rema}

\begin{rema}
	\label{Rk-Remainder-Thm4.1}
	The remainder terms $\|u \|_{L^2(\mathcal{L}_{T})}$ in \eqref{R-obs-t-x-H1-micro}, and $\| (u_0, u_1)\|_{L^2(\Omega) \times H^{-1}(\Omega)}$ in \eqref{R-obs-H1-micro} and \eqref{R-obs-H1-local}, can be weakened to
$\|u \|_{H^{-1}(\mathcal{L}_{T})}$ in \eqref{R-obs-t-x-H1-micro} and $\| (u_0, u_1)\|_{H^{-1}(\Omega)\times (H^2 \cap H^1_0(\Omega))'}$
respectively, using the same proof as in Theorem \ref{R-obs}. More generally, it is clear that these remainder terms can be replaced by norms in any Hilbert spaces of negative order, as long as they are appropriately adapted to the boundary conditions of the problem.

We also emphasize that this observation applies to the remainder term in the observability inequality \eqref{Distributed-Obs-Best-Geo-O-omega-H1}.

\end{rema}

\begin{rema}
	\label{Remark-micro-local-UC}In the proof of Theorem \ref{Theo-Obs-1}, we showed that the remainder terms in \eqref{R-obs-H1-local} and \eqref{R-obs-L2-local} can be removed through a simple analysis of the invisible set (see Lemma \ref{UC}), even in cases where the unique continuation property does not hold for all solutions of the wave equation \eqref{waveequation}. Further improvements along these lines are discussed in Remark \ref{Rk-Better-Local-Obs} and Corollary \ref{Cor-Obs-1-bis}.

It is natural to ask whether a similar strategy, as used in Lemma \ref{UC}, can be applied to remove the remainder terms in the microlocal estimates \eqref{R-obs-t-x-H1-micro}, \eqref{R-obs-H1-micro}, \eqref{R-obs-t-x-L2-micro}, and \eqref{R-obs-L2-micro}, at least for initial data microlocally supported in suitable regions. Unfortunately, this approach appears not to be effective in this context.

To illustrate the difficulty, consider the inequality \eqref{R-obs-t-x-H1-micro}. Let $\psi = \psi(t, x, D_t, D_x) \in \mathcal{A}^0$ be a pseudodifferential operator with \( \operatorname{supp} \psi \cap T^*_b \mathcal{L} \subset \mathscr{R}(\omega_T) \), and define the set
\[
\mathcal{N}_\psi = \left\{ u \in L^2_{\text{loc}}(\mathcal{L}) \ \middle| \
\begin{array}{l}
u \text{ solves } \eqref{waveequation}, \\
(I - \psi)u = 0 \text{ in } \mathcal{L}, \\
\partial_t u = 0 \text{ in } \omega_T
\end{array}
\right\}.
\]
By \eqref{R-obs-H1-micro}, any $u \in \mathcal{N}_\psi$ satisfies $u = \psi u \in H^1_{\text{loc}}(\mathcal{L})$, implying that $\mathcal{N}_\psi$ is compact and hence finite-dimensional.

Now consider $v = \partial_t u$ for $u \in \mathcal{N}_\psi$. Clearly, $v$ also solves the wave equation and satisfies $\partial_t v = 0$ in $\omega_T$. However, there is no guarantee that v belongs to $\mathcal{N}_\psi$, as we do not have $(I - \psi)v = 0$. In fact, since $\psi v = \psi (\partial_t u)$ and $v = \partial_t u$, their difference involves the commutator $[\psi, \partial_t] u$, which does not vanish in general. Therefore, the naive use of the operator $\partial_t$ does not yield an operator acting invariantly on $\mathcal{N}_\psi$.

We have not been able to further analyze the structure of the sets $\mathcal{N}_\psi$. Whether or not $\mathcal{N}_\psi$ is non-trivial remains an open problem.
\end{rema}

When, in addition, the unique continuation property holds for \eqref{waveequation}, i.e., when the uniqueness condition \eqref{Uniqueness-Condition} is satisfied, we can get the following  result:

\begin{theorem}[A second relaxed microlocal observability estimate] 
\label{Thm-distributed-observation-UC} 
		Let $\omega$ be a non-empty open set of $\Omega$ and $T >0$ such that 
	\begin{equation}
		\label{Uniqueness-Condition-2}
		T > 2 \sup_{\Omega} d(x, \omega).		
	\end{equation}
	\begin{enumerate}
		\item [(1a)] For every operator $\psi (t,x,D_t, D_{x}) \in \mathcal{A}^{0}$ with $\supp(\psi)\cap T^{*}_{b}\mathcal{L} \subset \mathscr{R}(\omega_{T})$,   there exists $C>0$ such that for every initial data $(u_{0},u_{1}) \in H^1_0(\Omega)\times L^{2}(\Omega)$,  the solution $u$ of  \eqref{waveequation} satisfies the observability estimate  
	\begin{equation}
		\label{Distributed-Obs-Best-H--1-t-x}
		\Vert \psi (t,x,D_t,D_{x}) u \Vert_{ H^{1}(\mathcal{L}_{T})} 
		\leq C\Vert \partial_t u \Vert_{L^{2}(\omega_{T} )} 
		+C \Vert  (I- \psi (t,x,D_t, D_{x}))u \Vert_{L^{2}(\mathcal{L}_{T})}.
	\end{equation}
	
		\item [(2a)] For every operator $\psi (x,D_{x}) \in \mathcal{B}^{0}$ with $\supp(\psi)\cap T^{*}_{b}\Omega \subset \mathscr{R}_0(\omega_{T})$,   there exists $C>0$ such that for every initial data $(u_{0},u_{1}) \in H^1_0(\Omega) \times L^{2}(\Omega)$,  the solution $u$ of  \eqref{waveequation} satisfies the observability estimate  
	\begin{equation}
		\label{Distributed-Obs-Best-H--1}
		\Vert \psi (x,D_{x})(u_{0},u_{1})\Vert_{ H^{1}_{0}(\Omega) \times L^{2}(\Omega)} 
		\leq C\Vert \partial_t u\Vert_{L^{2}(\omega_{T} )} +C \Vert  (I- \psi (x,D_{x}))(u_{0},u_{1})\Vert_{L^2(\Omega) \times  H^{-1}(\Omega) }.
	\end{equation}
	\end{enumerate}

	Similarly	
	\begin{enumerate}
		\item [(1b)] For every operator $\psi (t,x,D_t, D_{x}) \in \mathcal{A}^{0}$ with $\supp(\psi)\cap T^{*}_{b}\mathcal{L} \subset \mathscr{R}(\omega_{T})$,   there exists $C>0$ such that for every initial data $(u_{0},u_{1}) \in L^{2}(\Omega)\times H^{-1}(\Omega)$,  the solution $u$ of  \eqref{waveequation} satisfies the observability estimate  
	\begin{equation}
		\label{Distributed-Obs-Best-L--2-t-x}
		\Vert \psi (t,x,D_t,D_{x}) u \Vert_{ L^{2}(\L_{T})} 
		\leq C\Vert u \Vert_{L^{2}(\omega_{T} )} 
		+C \Vert  (I- \psi (t,x,D_t, D_{x}))u \Vert_{ H^{-1}(\L_{T})}.
	\end{equation}
	
		\item [(2b)] For every operator $\psi (x,D_{x}) \in \mathcal{B}^{0}$ with $\supp(\psi)\cap T^{*}_{b}\Omega \subset \mathscr{R}_0(\omega_{T})$,   there exists $C>0$ such that for every initial data $(u_{0},u_{1}) \in L^{2}(\Omega)\times H^{-1}(\Omega)$,  the solution $u$ of  \eqref{waveequation} satisfies the observability estimate  
	\begin{equation}
		\label{Distributed-Obs-Best-L--2}
		\Vert \psi (x,D_{x})(u_{0},u_{1})\Vert_{ L^{2}(\Omega) \times H^{-1}(\Omega)} 
		\leq C\Vert u\Vert_{L^{2}(\omega_{T} )} +C \Vert  (I- \psi (x,D_{x}))(u_{0},u_{1})\Vert_{ H^{-1}(\Omega) \times (H^2 \cap H^1_0)'(\Omega)}.
	\end{equation}
	\end{enumerate}
\end{theorem}

\begin{rema} Note that the main difference between  Theorem \ref{R-obs} and Theorem \ref{Thm-distributed-observation-UC} is that, in the latter, the compact remainder term is localized through the projection realized by the pseudodifferential operator $(I - \psi(t,x, D_t, D_x))$ or $(I- \psi (x,D_{x}))$, while, in the first theorem, the remainder involves the whole initial data. But for this to be done, we have assumed the condition \eqref{Uniqueness-Condition-2} guaranteeing the time-horizon is large enough to ensure that unique continuation holds. Whether the results of Theorem \ref{Thm-distributed-observation-UC} can be achieved from those in Theorem \ref{R-obs} without any additional further unique continuation assumption by means of a compactness-uniqueness argument is an interesting open problem, as we have discussed above in Remark \ref{Remark-micro-local-UC}.
\end{rema}

\begin{rema}
	Let us point out that, using Remark \ref{Rk-Remainder-Thm4.1}, we can weaken the remainder terms in the estimates of Theorem \ref{Thm-distributed-observation-UC}, replacing the terms $ \Vert  (I- \psi (t,x,D_t, D_{x}))u \Vert_{L^{2}(\mathcal{L}_T)}$ in \eqref{Distributed-Obs-Best-H--1-t-x} and $\Vert  (I- \psi (x,D_{x}))(u_{0},u_{1})\Vert_{L^2(\Omega) \times  H^{-1}(\Omega) }$ in \eqref{Distributed-Obs-Best-H--1} by the weaker terms  $\Vert  (I- \psi (t,x,D_t, D_{x}))u \Vert_{H^{-1}(\mathcal{L}_T)}$ in \eqref{Distributed-Obs-Best-H--1-t-x} and $\Vert  (I- \psi (x,D_{x}))(u_{0},u_{1})\Vert_{ H^{-1}(\Omega) \times (H^2 \cap H^1_0(\Omega))'}$ in \eqref{Distributed-Obs-Best-H--1}.
\end{rema}

\begin{rema}[Examples]\label{+examples}The typical example in which Theorem 
\ref{Thm-distributed-observation-UC} applies is for instance when $\Omega$ is the unit ball, $A = Id$ (i.e.,  the standard constant coefficients wave equation),  and $\omega$ is the ball of radius $1/2$. In such case, the GCC is not satisfied in any time, due to the whispering gallery phenomenon, i.e., the existence of rays localized in a neighborhood of the boundary. However, as soon as $T > 1$, unique continuation holds and the above theorem applies. 

Another example corresponds to the case in which $\Omega$ is the unit square, $A = Id$, and $\omega$ is an $\varepsilon$($>0$)-neighborhood of the whole boundary. Note that,  even if the square is not a smooth bounded domain and its boundary has tangencies of infinite order, this is not an impediment for our results to apply  since the observation is made on a neighborhood of the whole boundary and, therefore, boundary phenomena are irrelevant. In this case, the GCC holds as soon as $T > \sqrt{2} (1 - 2 \varepsilon)$ while unique continuation holds as soon as $T > (1 - 2 \varepsilon)$. Therefore, when $T$ belongs to the intermediate interval $((1 - 2\varepsilon), \sqrt{2} (1 - 2 \varepsilon))$, Theorem \ref{Thm-distributed-observation-UC} applies, but the classical  observability inequality based on GCC does not hold. 

The same occurs in most domains $\Omega$ since a $\varepsilon$-neighborhood of the boundary always guarantees  GCC  when the time-horizon is long enough, but, normally, the unique continuation property holds in shorter times.
\end{rema}

We conclude this section with the following result, which goes a step further by addressing the case where global unique continuation fails. In such situations, it becomes necessary to assume that either the initial position $u_0$ or the initial velocity $u_1$ vanishes.

\begin{theorem}[When unique continuation does not hold and $u_0$ or $u_1$ vanishes]
	\label{Thm-Obs-Without-UC}
	Let $ \Omega$ be a smooth (possibly unbounded) domain, $\omega$  a non-empty bounded open set of $\Omega$ and $T >0$.

	Then all the items (1a), (2a), (1b) and (2b) in Theorem \ref{Thm-distributed-observation-UC} hold for solutions $u$ of the wave equation \eqref{waveequation} corresponding to initial data $(u_0, u_1)$ with either $u_0 = 0$ or $u_1 = 0$.
\end{theorem}

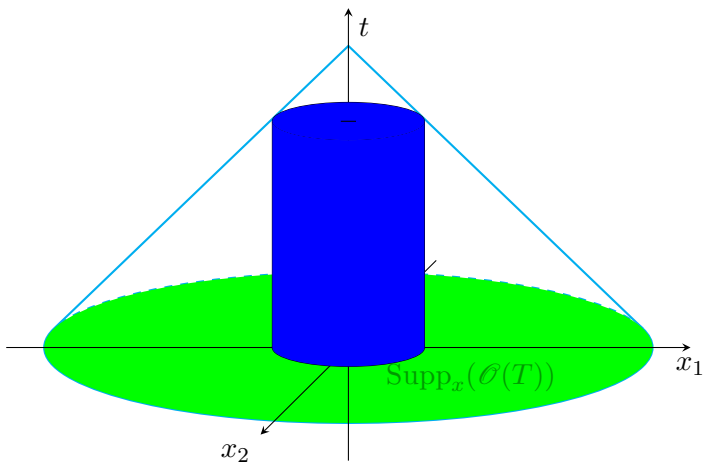
\begin{figure}
\label{Fig-R-d}
\begin{tikzpicture}[>=stealth,join=round]
\def\a{4}   
\def\b{1}  
\def\h{4}   
\def\d{2}   
\pgfmathsetmacro{\t}{asin(\b/\h)}  

\begin{scope}[cyan,thick]

\draw[dashed]
(\t:{\a} and {\b}) arc(\t:180-\t:{\a} and {\b});
\draw
(\t:{\a} and {\b})--(0,\h)--(180-\t:{\a} and {\b})
arc(180-\t:360+\t:{\a} and {\b});

\end{scope}

    \fill[green, opacity=0.3] (4,0)  arc (0:360:4 and 1) ;
    \node[green!60!black] at (2, 0, 1) {$\supp_x(\mathcal{O}(\omega_T))$};

    \draw[->] (-4.5,0,0) -- (4.5,0,0) node[below] {$x_1$};
    \draw[->] (0,-1.5,0) -- (0,4.5,0) node[below right] {$t$};
    \draw[->] (0,0,-3) -- (0,0,3) node[below left] {$x_2$};

\draw (-1,0) arc (180:360:1 and 0.25);
\draw[dashed] (-1,0) arc (180:360:1 and -0.25);
\draw (-1,0) -- (-1,3);
\draw (-1,3) arc (180:360:1 and 0.25);
\draw (-1,3) arc (180:360:1 and -0.25);
\draw (1,3) -- (1,0);  
\fill [blue, opacity=0.2] (-1,0) -- (-1,3) arc (180:360:1 and 0.25) -- (1,0) arc (0:180:1 and -0.25);
\fill [blue, opacity=0.1] (1,3) arc (0:360:1 and -0.25);
 \node[blue] at (0.5,3,0) {$T$} ; 
\draw (-0.1,3,0) -- (0.1,3,0) ;

\end{tikzpicture}
 \caption{When $\Omega = \R^d$, $ \omega $ is the unit ball and $\Delta_A$ is the flat Laplacian, the set $\mathcal{O}(\omega_T)$, when projected in the physical space is supported in the ball of size $T$. In fact, in this case, as the bicharacteristics in this setting correspond to straight lines, $$\mathcal{O}(\omega_T) = \{(x, \xi) \in T_{b}^{*}\Omega\backslash 0 \text{ s.t. } \exists t \in [0,T] \text{ satisfying } x + \xi t \in \omega \}.$$ }
\end{figure}

The typical example in which Theorem \ref{Thm-Obs-Without-UC} applies is when $\Omega = \R^d$ and $ \omega $ is the unit ball, see Figure \ref{Fig-R-d}.

\subsection{Proofs}

\begin{proof}[Proof of Theorem \ref{R-obs}]
All the estimates of Theorem \ref{R-obs} can be proved through a direct application of the closed graph theorem and Lemma \ref{regularity-general}. Below, we present the proof of the estimate \eqref{R-obs-H1-micro}, the other proofs  being completely similar and left to the reader.

We will closely follow the proof of estimate \eqref{R-obs-L1} of Corollary \ref{relaxed-obs}. Consider  the following Hilbert space 
\begin{equation*}
	E'=\Big\{ (u_{0},u_{1}) \in L^{2}(\Omega)\times H^{-1}(\Omega), \,\, \text{ and } \d_{t}u \in L^{2}(\omega_{T})\Big\}
\end{equation*}
equipped with the norm 
$$\Vert (u_{0},u_{1})\Vert_{E'}^{2}=\Vert (u_{0},u_{1})\Vert_{L^{2}(\Omega)\times H^{-1}(\Omega)}^{2} + \Vert \d_{t}u\Vert_{L^{2}(\omega_{T})}^{2},$$ 
and the  energy space $F = H_{0}^{1}(\Omega)\times L^{2}(\Omega)$
equipped with its natural norm. According to Lemma \ref{regularity-general}, item 2, for every operator $\psi (x,D_{x}) \in \mathcal{B}^{0}$ such that  $\supp(\psi)\cap T^{*}_{b}\Omega \subset \mathcal{R}_{0}(\omega_{T})$,  the  map 
\begin{equation*}
	E' \longrightarrow F = H_{0}^{1}(\Omega)\times L^{2}(\Omega) ,  \quad (u_{0},u_{1}) \mapsto \psi (x,D_{x})(u_{0},u_{1})
\end{equation*}
is well defined. Consequently,  the closed graph theorem  yields its continuity and estimate \eqref{R-obs-H1-micro}.
\end{proof}

\begin{proof}[Proof of Theorem \ref{Thm-distributed-observation-UC}]
	Here again, all the items of Theorem \ref{Thm-distributed-observation-UC} can be proved similarly using a classical compactness-uniqueness argument. Below, we only present the proof of the estimate \eqref{Distributed-Obs-Best-H--1}, as the other ones follow exactly the same path.
	
	Let $\psi (x,D_{x}) \in \mathcal{B}^{0}$ with $\supp(\psi)\cap T^{*}_{b}\Omega \subset \mathscr{R}_0(\omega_{T})$.
In view of \eqref{R-obs-H1-micro} it is sufficient to show the existence of a constant $C>0$ such that 
\begin{equation}
\label{Main-Step-Compactness-Uniqueness}
\Vert   \psi (x,D_{x})(u_{0},u_{1})\Vert_{L^{2}(\Omega)\times H^{-1}(\Omega)} \le C \left( \Vert \d_{t}u\Vert_{L^{2}(\omega_{T} )} + \Vert  (I- \psi (x,D_{x}))(u_{0},u_{1})\Vert_{L^{2}(\Omega)\times H^{-1}(\Omega)}\right),
\end{equation}
for all solution $(u_0, u_1) \in H^1_0 (\Omega) \times L^2(\Omega)$.

We argue by contradiction. If that were not the case it would exist a sequence $(u_{0,k}, u_{1,k})_{k \in \N} \subset H^1_0 (\Omega) \times L^2(\Omega)$ such that 
	\begin{equation}
		\label{Contradiction-Assumptions}
		 \Vert \psi (x,D_{x})(u_{0,k},u_{1,k})\Vert_{L^{2}(\Omega)\times H^{-1}(\Omega)} = 1, 
	\end{equation}
	\begin{equation}
		 \lim_{k \to \infty} \left( \Vert \d_{t}u_k\Vert_{L^{2}(\omega_{T} )} + \Vert  (I- \psi (x,D_{x}))(u_{0,k},u_{1,k})\Vert_{L^{2}(\Omega)\times H^{-1}(\Omega)}\right) = 0.
	\end{equation}
	Accordingly, $(u_{0,k},u_{1,k})_{k \in \N}$ is bounded in $L^2(\Omega) \times H^{-1}(\Omega)$, and up to the extraction of a subsequence still denoted the same, weakly converges to some $(u_0, u_1)$, and from the above condition, we also have that $(I - \psi(x,D_x)) (u_0, u_1) = 0$, that is $(u_0,u_1) = \psi(x, D_x)(u_0, u_1)$.
	
	Then, in view of \eqref{R-obs-H1-micro}, $\psi(x, D_x)(u_{0,k},u_{1,k})$ is bounded in $H^1_0(\Omega) \times L^2(\Omega)$, so it weakly converges to $\psi(x, D_x) (u_0, u_1) = (u_0, u_1) $ in $H^1_0(\Omega) \times L^2(\Omega)$, entailing in particular that $(u_0, u_1) \in H^1_0(\Omega) \times L^2(\Omega)$, and corresponds to a solution $u$ of \eqref{waveequation} such that $\partial_t u = 0$ in $\omega_{T}$.

	By unique continuation we deduce that the limit $u\equiv 0$ and therefore $(u_0, u_1) \equiv (0,0)$. But then the sequence $(\psi(x,D) (u_{0,k},u_{1,k}))_{k \in \N}$ strongly converges to $(0,0)$ in  $L^2(\Omega) \times H^{-1}(\Omega)$. This contradicts \eqref{Contradiction-Assumptions} and concludes the proof.
\end{proof}

\begin{proof}[Proof of Theorem \ref{Thm-Obs-Without-UC}]
	The proof of Theorem \ref{Thm-Obs-Without-UC} follows the one of Theorem \ref{Thm-distributed-observation-UC}, the only difference being the unique continuation property we shall rely on, which is the following one: if $u$ is a solution of the wave equation \eqref{waveequation} corresponding to an initial datum satisfying $u_0 = 0$ (respectively $u_1 = 0$) such that $\partial_t u = 0 $ in  $(0,T) \times \omega$, then $u$ vanishes identically in $\{ (t,x) \in (-T, T) \times \Omega, \, d(x,\omega) + |t| \leq T \}$, and thus $u_0$ (respectively $u_1$) vanishes in the set $O^{T} = \{ x \in \overline{\Omega}, \, d(x, \omega) < T\}$.
	
	Indeed, if $u_0$ (respectively $u_1$) vanishes, the function $u$ extended in a odd  (respectively even) manner is a solution of the wave equation \eqref{waveequation} on $(-T, T) \times \Omega$. We can then use the classical unique continuation result for the wave equation \cite{Tataru95} which asserts that, if $\partial_t u = 0$ in $(-T,T) \times \omega$ for a solution $u$ of \eqref{waveequation} on $(-T, T) \times \Omega$, then $u$ vanishes in the set $\{ (t,x) \in (-T, T) \times \Omega, \, d(x,\omega) + |t| \leq T \}$.
	
	Let us now explain how it can be used to prove for instance (again, all the other statements in Theorem \ref{Thm-Obs-Without-UC} can be proved similarly) that for every operator $\psi (x,D_{x}) \in \mathcal{B}^{0}$ with $\supp(\psi)\cap T^{*}_{b}\Omega \subset \mathscr{R}_0(\omega_{T})$,  there exists $C>0$ such that for  any initial data $(u_0, u_1)$ with $u_0 = 0$ and $u_1 \in L^2(\Omega)$, the solution of  \eqref{waveequation} satisfies the observability estimate  
	\begin{equation}
		\label{Distributed-Obs-Best-without-UC-1}
	\Vert \psi (x,D_{x})u_{1}\Vert_{H^{-1}(\Omega)} 
		\leq C\Vert \d_{t}u\Vert_{L^{2}(\omega_{T} )} +C \Vert  (I- \psi (x,D_{x}))u_{1}\Vert_{H^{-1}(\Omega)}.
	\end{equation}

We mimic the proof of the estimate \eqref{Distributed-Obs-Best-H--1} of Theorem \ref{Thm-distributed-observation-UC}, and use a compactness uniqueness argument to prove that \eqref{Main-Step-Compactness-Uniqueness} holds for any initial data $(u_0, u_1)$ with $u_0 = 0$ and $u_1 \in L^2(\Omega)$. 

By contradiction and following the proof of the estimate \eqref{Distributed-Obs-Best-H--1} of Theorem \ref{Thm-distributed-observation-UC}, we get a sequence $u_{1,k}$ such that  $   ((I- \psi (x,D_{x}))u_{1,k})_{k \in \N} $ goes to $0$ in $H^{-1}(\Omega)$, $\psi(x,D_x) u_{1,k}$ is of unit norm in $H^{-1}(\Omega)$, and such that the corresponding solutions $u_k$ of \eqref{waveequation} satisfies that $(\partial_t u_k)_{k \in \N}$ goes to $0$ in $L^2(\omega_{T})$. Consequently, up to a subsequence, we get $u_1$ such that $(u_{1,k})$ converges weakly to $u_1$ in $H^{-1}(\Omega)$, and $(I- \psi (x,D_{x}))u_{1}= 0$, and such that the corresponding solution of \eqref{waveequation} with initial data $(0, u_1)$ satisfies $\partial_t u = 0$ in $(0,T) \times \omega$. By the above uniqueness result, we thus get that $u_1 = 0$ in the set $O^{T} = \{ x \in \overline{\Omega}, \, d(x, \omega) < T\}$. Finally, since $\psi = \psi (x,D_{x}) \in \mathcal{B}^{0}$ satisfies $\supp(\psi)\cap T^{*}_{b}\Omega \subset \mathscr{R}_0(\omega_{T})$ and the $x$-projection of $ \mathscr{R}_0(\omega_{T})$ is included in $O^{T}$, $u_1 = \psi(x, D_x) u_1$ implies that $u_1$ is supported in $O^{T}$. Therefore, $u_1$ vanishes everywhere. Now, using \eqref{R-obs-H1-micro}, $\psi(x, D_x)(u_{1,k})$ is bounded in $ L^2(\Omega)$, so we also obtain by compactness that $(\psi(x,D_x) u_{1,k})_{k \in \N}$ strongly converges to $0$ in $H^{-1}(\Omega)$, thus getting a contradiction.
\end{proof}

 \section{Control theoretical consequences}\label{Sec-Control-without-GCC}

Each of the new observability results we have presented have their counterpart at the control level. This can be seen systematically by the employment of the duality arguments as  in  \cite{Lions,LionsSIAM88}.

Note however that, duality  transfers the observability of the adjoint backward wave equation into the control of the forward wave process. Thus, attention has to be paid to rewriting the needed microlocal assumptions reversing the sense of time. This is a purely technical minor aspect since we are dealing with  time-independent variable coefficients and the geometry of the relevant pairs $(\omega, \overline{\mathscr{O}})$ is independent of the sense of time. Waves with time-dependent coefficients pose new technical difficulties,  as we will discuss in the last section.

\subsection{Controllable $(\omega, \overline{\mathscr{O}})$ pairs}

The following result is the counterpart of Theorem \ref{Theo-Obs-1} from the control point of view:

\begin{theorem}\label{Thm-Obs-local}
	
	Within the setting of Theorem \ref{Theo-Obs-1},   for  every data $(y_{0}^T, y_{1}^T) \in H_{0}^{1}(\Omega)\times L^{2}(\Omega)$, there exists a control $v \in L^2(0,T; L^2(\omega))$ such that the  solution $y$ of  
	\begin{equation}
	\left\{ 
	\begin{array}{ll}
		\partial_{t}^{2}y -\Delta_{A}y =v 1_{\omega}\quad& \text{in } (0,T) \times \Omega ,
		\\
		y(t, \cdot)=0\quad &\text{on } (0,T) \times \partial \Omega, 
		\\
		(y(0, \cdot ),\partial _{t} y(0,\cdot))= (0, 0) & \text{in } \Omega.
\end{array}%
\right.  \label{waveequation-control}
	\end{equation}
	satisfies
	\begin{equation}
		\label{Control-Req-y-Local}
		y(T, \cdot) = y_0^T  \quad \text{ and }\quad \partial_t y(T, \cdot) =  y_1^T\quad  \text{ in }\,\, \overline{\mathscr{O}}.
	\end{equation}
	Furthermore, there exists $C>0$ such that
	\begin{equation}
		\vert| v \vert|_{L^2(\omega_{T})}
		\leq  
		C \| (y_0^T, y_1^T)\|_{H^1_0(\Omega) \times L^2(\Omega)}.
	\end{equation}
\end{theorem}

\begin{proof}
	Since the set $\mathcal{O}(\omega_{T})$ is open and $\overline{\mathscr{O}} \subset \mathcal{O}(\omega_{T}) $ by assumption, there exists an open set $\mathscr{O}_1$ such that $\overline{\mathscr{O}} \subset \mathscr{O}_1$ and $\overline{\mathscr{O}_1} \subset \mathcal{O}(T)$. We then take $\chi \in \mathscr{C}^\infty_c (\mathscr{O}_1)$ which equals to $1$ in $\overline{\mathscr{O}}$. 

	Applying the observability inequality \eqref{obs-theo-L2} of Theorem \ref{Theo-Obs-1} on $\mathscr{O}_1$, we get that for any initial data $(u_0, u_1) \in L^2(\Omega) \times H^{-1}(\Omega)$ supported in $\overline{\mathscr{O}_1}$, 
	\begin{equation}
		\label{Obs-O1}
		\| (u_0, u_1) \|_{L^2(\Omega) \times H^{-1}(\Omega) } 
		\leq 
		C \| u \|_{L^2(\omega_{T})}, 
	\end{equation}
	where $u$ is the corresponding solution of \eqref{waveequation}.
	
	By time reversal ($t \mapsto T-t$), for any initial data $(u_0, u_1) \in L^2(\Omega) \times H^{-1}(\Omega)$ supported in $\overline{\mathscr{O}_1}$, the solution of 
	\begin{equation}
	\left\{ 
	\begin{array}{ll}
	\partial _{t}^{2}u-\Delta_{A}u=0\quad &\text{in } \,\,\L 
	\\ 
	u(t,.)=0\quad &\text{on } \,\,\d\L 
	\\
	(u(T,\cdot),\partial _{t}u(T,\cdot))=(u_{0}^T,u_{1}^T), & 
	\end{array}%
	\right.  
	\label{waveequation-backward}
\end{equation}
satisfies
\begin{equation}
		\label{Obs-O1-T}
		\| (u_0^T, u_1^T) \|_{L^2(\Omega)  \times H^{-1}(\Omega) } 
		\leq 
		C \| u \|_{L^2(\omega_{T})}, 
	\end{equation}

	We then introduce the set 
	$$
		X = \{(u_0, u_1) \in L^2(\Omega) \times H^{-1}(\Omega) \text{ supported in $\overline{\mathscr{O}}_1$} \}, 
	$$
	which is obviously closed for the $L^2 \times H^{-1}$ topology. 
	
	Take $(y_0^T, y_1^T) \in H^1_0(\Omega) \times L^2(\Omega)$, and introduce the functional $J$ defined for $(u_0^T, u_1^T) \in X$ by 
	\begin{equation}
		J (u_0^T, u_1^T) 
		= \frac{1}{2} \int_0^{T} \int_\omega | u |^2 \dx \dt
		- 
		\int_\Omega \chi u_0^T  y_1^T \dx
		+
		\langle \chi u_1^T,  y_0^T \rangle_{H^{-1}(\Omega) , H^1_0(\Omega)}, 
	\end{equation}
	where $u$ is the corresponding solution of \eqref{waveequation-backward}.
	
	Here and it what follows $\langle \cdot,  \cdot \rangle_{H^{-1}(\Omega) , H^1_0(\Omega)}, $ stands for the duality pairing between $H^{-1}(\Omega) $ and $H^1_0(\Omega)$.

	It is obvious from the estimate \eqref{Obs-O1} that $J$ is continuous, strictly convex and coercive on $X$. Therefore, there exists a minimizer $(U_0, U_1) \in X$ of $J$ such that 
	$$
		\| U \|_{L^2(\omega_{T})} \leq C \| (y_0, y_1) \|_{H^1_0(\Omega)  \times L^2(\Omega) }.
	$$
	The Euler Lagrange equation then gives that for all $(u_0, u_1) \in X$, 
	$$
		0 = \int_0^{T} \int_\omega U u  \dx \dt
		-
		\int_\Omega \chi u_0  y_1 \dx
		+
		\langle \chi u_1,  y_0 \rangle_{H^{-1}(\Omega) , H^1_0(\Omega) }. 
	$$
	Since the solution $y$ of \eqref{waveequation-control} corresponding to a control function $v \in L^2(0,T; L^2(\omega))$ satisfies that for all $(u_0, u_1) \in X$, 
	$$
		0 = \int_0^{T} \int_\omega v u  \dx \dt
		-
		\int_\Omega  u_0  \partial_t y(T) \dx
		+
		\langle u_1,  y(T,\cdot) \rangle_{H^{-1}(\Omega) , H^1_0(\Omega) }, 
	$$
	by setting
	$$
		v = U 1_\omega, 
	$$
	we observe that the corresponding solution $y$ of \eqref{waveequation-control} satisfies
	$$
		y(T, \cdot ) = \chi y_0 \text{ and } \partial_t y(T, \cdot) = \chi y_1 \,\, \text{ in } \mathscr{O}_1. 
	$$
	This concludes the proof of Theorem \ref{Thm-Obs-local}.
\end{proof} 	

\begin{rema}
	\label{Remark-Using-Cor}
	Starting from Corollary \ref{Cor-Obs-1-bis}, we can improve the result of Theorem \ref{Thm-Obs-local} as follows. For $\chi \in \mathscr{C}^\infty_c(\mathcal{O}(\omega_{T}))$,   for  every data $(y_{0}^T, y_{1}^T) \in H_{0}^{1}(\Omega)\times L^{2}(\Omega)$, there exists a control $v \in L^2(0,T; L^2(\omega))$ such that the solution $y$ of  \eqref{waveequation-control} satisfies
	\begin{equation}
		\label{Control-Req-y-Local-chi}
		y(T, \cdot) = y_0^T  \quad \text{ and }\quad \partial_t y(T, \cdot) =  y_1^T\quad  \text{ in }\,\, \{\chi = 1\}.
	\end{equation}
	and there exists $C>0$ such that
	\begin{equation}
		\vert| v \vert|_{L^2(\omega_{T})}
		+
		\| (y(T), \partial_t y(T)) - \chi (y_0^T, y_1^T)\|_{H^2 \cap H^1_0(\Omega) \times H^1_0(\Omega)} 
		\leq  
		C \| (y_0^T, y_1^T)\|_{H^1_0(\Omega) \times L^2(\Omega)}.
	\end{equation}
	Note that, since $v$ belongs to $L^2(\omega_{T})$, we should rather expect the solution $y$ of \eqref{waveequation-control} to be in $\mathscr{C}^0([0,T]; H^1_0(\Omega)) \cap \mathscr{C}^1 ([0,T]; L^2(\Omega))$. In other words, such improvement means that we can construct a control process that controls exactly the solution $y$ at time $T$ on $\overline{\mathcal{O}}$ and do not create $H^1$ singularities outside of the support of $1 - \chi$.
	
	In order to prove such result, simply replace the functional $J$ above by $J_\chi$ defined by 
	\begin{multline}
		J_{\chi} (u_0^T, u_1^T) 
		= \frac{1}{2} \int_0^{T} \int_\omega | u |^2 \dx \dt
		+ \frac{1}{2} \| (1- \chi) (u_0^T, u_1^T)\|_{H^{-1} \times H^{-2}}^2 
		\\
		- 
		\int_\Omega \chi u_0^T  y_1^T \dx
		+
		\langle \chi u_1^T,  y_0^T \rangle_{H^{-1}(\Omega) , H^1_0(\Omega)}, 
	\end{multline}
	for $(u_0^T, u_1^T) \in L^2(\Omega) \times H^{-1} (\Omega)$,  
	where $u$ is the corresponding solution of \eqref{waveequation-backward}, and $H^{-2}$ is a short notation for $(H^2 \cap H^1_0(\Omega))'$.

	The observability estimate \eqref{Distributed-Obs-Best-Geo-O-omega-L2} easily provides the coercivity and strict convexity of the functional $J_\chi$ on the space $X_{obs} = \overline{L^2(\Omega) \times H^{-1} (\Omega)}^{\|\cdot \|_{obs}}$, where the norm $\| \cdot \|_{obs}$ is given by 
	$$
		\|  (u_0^T, u_1^T) \|_{obs}^2 = \int_0^{T} \int_\omega | u |^2 \dx \dt
		+\| (1- \chi) (u_0^T, u_1^T)\|_{H^{-1} \times H^{-2}}^2. 
	$$
	There is therefore a unique minimizer $(U_0^T, U_1^T) \in X_{obs}$ of $J_\chi$, which satisfies 
	$$
		\|  (U_0^T, U_1^T) \|_{obs} \leq C \| (y_0^T, y_1^T) \|_{H^1_0(\Omega) \times L^2(\Omega)}. 
	$$
	Following the above proof, one then easily derives that 
	\begin{align*}
		y(T) &= \chi y_0 + (1- \chi) (-\Delta)^{-2} (1- \chi) U_1^T, 
		\\
		\partial_t y(T) & = \chi y_1 - (1- \chi) (-\Delta)^{-1} (1 - \chi) U_0^T, 
	\end{align*}
	from which we directly conclude the proof of the above statement.
\end{rema}

\begin{rema}
\label{Remark-Exact-Control+GlobalApprox}
When, in addition to the geometric conditions of Theorem \ref{Theo-Obs-1}, the time horizon $T$ is long enough so that unique continuation holds i.e.,condition \eqref{Uniqueness-Condition-2}, the control result above can be improved to guarantee the simultaneous approximate controllability and the control of the projections as in \eqref{Control-Req-y-Local}. More precisely,  for all $\varepsilon >0$ there exists a control $v_\varepsilon$ such that the solution satisfies both \eqref{Control-Req-y-Local} and 
$$
\| y(T, \cdot ) - y_0^T\|_{H^1_0(\Omega)} +  \| \partial_t y(T, \cdot) -  y_1^T \|_{L^2(\Omega)} \le \varepsilon. 
$$
To prove it, it suffices to minimise the functional $J_\varepsilon$ defined by 
	\begin{multline}
		J_\varepsilon (u_0^T, u_1^T) =
		\frac{1}{2} \int_0^T \int_\omega | u |^2 \dx \dt + \varepsilon
		 \| ((1-\chi ) u_0^T, (1-\chi )u_1^T) \|_{L^2(\Omega)  \times H^{-1}(\Omega) }
		\\
		- 
		\int_\Omega u_0^T  y_1^T \dx
		+
		\langle  u_1^T,  y_0^T \rangle_{H^{-1}(\Omega) , H^1_0(\Omega)}, 
	\end{multline}
	on $L^2(\Omega) \times H^{-1}(\Omega)$, following the arguments in Section 2 of \cite{zuazua1997finite}, to prove the coercivity of the functional $J_\varepsilon$ in $L^2(\Omega) \times H^{-1}(\Omega)$ and then writing the Euler Lagrange equation for the minimizer to deduce the control. Note, however, that this approach does not provide a quantitative estimate for the cost of controllability in this setting, i.e., on the norm of the control in terms of $\varepsilon$.
\end{rema}

\subsection{Pseudodifferential control when unique continuation holds} 
Rather than presenting all the control results that can be derived by duality from the observability estimates in Section \ref{sec-further results}, we focus below on a representative control result of microlocal nature, which serves as the counterpart to item (2) of Theorem \ref{Thm-distributed-observation-UC}.

In order to do so, for $T>0$ we introduce the set 
\begin{multline}\label{GCC-1-reversal}
\widetilde{\mathscr{R}_0(\omega_{T} ) }= \Big\{(x,\xi) \in T_{b}^{*}\Omega\backslash 0, \text{ such that  bicharacteristics  $\gamma_\rho$}
\\ \text{ issued from $(x,\xi)$ at time $T$ satisfy } \gamma_{\rho}(\R) \cap T^{*}(\omega_{T}) \neq \emptyset \Big\}.
 \end{multline}
 
Note that the set $\widetilde{\mathscr{R}_0(\omega_{T} ) }$ differs from $\mathscr{R}_0(\omega_{T} )$ in that it considers bicharacteristics originating from $(x, \xi)$ at time $T$,  rather than at the initial time. By a simple time-reversal argument  (i.e., the change of variable $t \mapsto T-t$), this is precisely the relevant set when the goal is to obtain information about $(u, \partial_t u)$ at time $t=T$ rather than at the initial time  $t=0$, for solutions of \eqref{waveequation}.

\begin{theorem}[{Pseudodifferential control}]
	\label{Thm-Control-With-UC}
	Assuming the uniqueness condition \eqref{Uniqueness-Condition}, for every operator $ \psi (x,D_{x}) \in \mathcal{B}^{0}$ with $\supp(\psi)\cap T^{*}_{b}\Omega \subset \widetilde{\mathscr{R}_0(\omega_{T})}$, there exists $C>0$ such that for  any initial data $(y_{0}^T, y_{1}^T) \in H_{0}^{1}(\Omega)\times L^{2}(\Omega)$, there exists $v \in L^2(0,T; L^2(\omega))$ such that the control $v$ and the corresponding solution $y$ of  \eqref{waveequation-control}
	 satisfies the following estimates: 
	\begin{multline}
		\label{Control-Req-y}
		\| (y(T, \cdot), \partial_t y(T, \cdot)) 
			-
		\psi(x, D_x)^\star (y_0^T, y_1^T) \|_{H^2 \cap H^1_0(\Omega) \times H^1_0(\Omega)}
		+ 
		\| v \|_{ L^2(\omega_T)}
		\\
		\leq  
		C \| (y_0^T, y_1^T)\|_{H^1_0(\Omega) \times L^2(\Omega)}.
	\end{multline}
\end{theorem} 

\begin{rema}
	Let us briefly comment the control requirement \eqref{Control-Req-y}. Here, let us emphasize that the target state $(y_0^T, y_1^T)$ belongs to $H^1_0(\Omega) \times L^2(\Omega)$ and the control function $v$ belongs to $L^2(\omega_T)$, so that the solution $(y, \partial_t y)$ of \eqref{waveequation-control} belongs to  $\mathscr{C}^0([0,T]; H^1_0(\Omega)) \times \mathscr{C}^1([0,T];L^2(\Omega))$. The relevant information of \eqref{Control-Req-y} is thus that we can choose a control function $v$ such that $ (y(T, \cdot), \partial_t y(T, \cdot)) -
		\psi(x, D_x)^\star (y_0^T, y_1^T)$ belongs to $H^2 \cap H^1_0(\Omega) \times H^1_0(\Omega)$, that is such that the $H^2 \cap H^1_0(\Omega) \times H^1_0(\Omega)$ singularities of $(y(T), \partial_t y(T))$ coincide with the ones of $\psi(x, D_x)^\star (y_0^T, y_1^T)$.
\end{rema}

\begin{proof} 
	Let $(y_0^T, y_1^T) \in H^1_0(\Omega) \times L^2(\Omega)$. We introduce the functional 
	\begin{multline}
		J (u_0^T, u_1^T) = \frac{1}{2} \int_0^T \int_\omega |\partial_t u |^2 \dx \dt
		+ 
		\frac{1}{2} \| (I - \psi(x, D_x))(u_0^T, u_1^T) \|_{L^2\times H^{-1}}^2 
		\\
		- 
		\int_\Omega A \nabla u_0^T \cdot \nabla \psi(x, D_x)^\star  y_0^T\dx
		- 
		\int_\Omega u_1^T \psi(x, D_x)^\star  y_1^T \dx, 
	\end{multline}
	defined for $(u_0, u_1) \in H^1_0(\Omega) \times L^2(\Omega)$, where $u$ is the corresponding solution of \eqref{waveequation-backward}.
	
	Here, to be precise, we define the $H^{-1}(\Omega)$-norm by the formula
	$$
		\| f \|_{H^{-1}(\Omega)}^2 = \int_\Omega A \nabla (-\Delta_A)^{-1} f \cdot \nabla (-\Delta_A)^{-1} f \dx, 
	$$
	where $-\Delta_A$ is the operator $-\div(A \nabla \cdot )$ in $\Omega$ with domain $ H^1_0(\Omega)$ on $H^{-1}(\Omega)$.
	
	From \eqref{Distributed-Obs-Best-H--1}, it is clear that the quantity
	$$
		\|(u_0^T, u_1^T)  \|_{obs}^2 
		=
		\int_0^T \int_\omega |\partial_t u |^2 \dx \dt
		+ 
		\| (I - \psi(x, D_x))(u_0^T, u_1^T) \|_{L^2\times H^{-1}}^2 
	$$
	defines a norm on $H^1_0(\Omega) \times L^2(\Omega)$, and we consider the closure $X$ of $H^1_0(\Omega) \times L^2(\Omega)$ with respect to this norm.  Note that we easily have 
	$$
		\| (u_0^T, u_1^T) \|_{L^2 \times H^{-1}} \leq C \| (u_0^T, u_1^T)\|_{obs}.
	$$

	We then check that the linear maps 
	$$
		(u_0^T,u_1^T) \mapsto \int_\Omega A \nabla u_0^T \cdot \nabla \psi(x, D_x)^\star  y_0^T \dx
		\quad \text{ and } \quad 
		(u_0^T,u_1^T) \mapsto \int_\Omega u_1 \psi(x, D_x)^\star  y_1 \dx
	$$
	are continuous with respect to the norm $\| \cdot \|_{obs}$: Indeed, 
	$$
		\left|  \int_\Omega u_1^T \psi(x, D_x)^\star  y_1^T \dx\right|
		= 
		\left|  \int_\Omega  \psi(x, D_x)  u_1^T y_1^T \dx \right|
		\leq 
		\| \psi(x, D_x) u_1^T \|_{L^2} \| y_1^T\|_{L^2}
		\leq 
		C \| (u_0^T, u_1^T)\|_{obs} \| y_1^T\|_{L^2}, 
	$$
	and
	\begin{align*}
		&\left|  
		\int_\Omega A \nabla u_0^T \cdot \nabla \psi^\star(x, D_x)  y_0^T \dx
		\right|
		 \leq 
		\left|  
		\int_\Omega A \nabla \psi(x, D_x) u_0^T \cdot \nabla y_0^T\dx
		\right|
		+ 
		\left| \langle  u_0^T, [\psi(x,D_x)^\star, \div (A \nabla \cdot) ]y_0^T \rangle \right| 
		\\ 
		&\qquad  \leq 
		C \| \psi(x, D_x) u_0^T \|_{H^1_0} \| y_0^T \|_{H^1_0}
		+ 
		C \| u_0^T \|_{L^2} \| y_0^T \|_{H^1_0}
		\leq C \| (u_0^T, u_1^T)\|_{obs} \| y_0^T\|_{H^1_0}.		
	\end{align*}

	Accordingly, the functional $J$ can be extended uniquely as a continuous coercive functional on $X$, and it has a unique minimizer $(U_0^T, U_1^T) \in X$, which satisfies
	$$
		\| (U_0^T, U_1^T) \|_{obs} \leq C \| (y_0^T, y_1^T) \|_{H^1_0 \times L^2}.
	$$
	
	The Euler-Lagrange equation satisfied by $(U_0,^T U_1^T)$ then gives that for all $(u_0^T, u_1^T) \in H^1_0(\Omega) \times L^2(\Omega)$, 
	\begin{multline*}
		0=  \int_0^T \int_\omega \partial_t U \partial_t u \dx \dt
		+ 
		\langle (I - \psi(x, D_x))(U_0^T, U_1^T), (I - \psi(x, D_x))(u_0^T, u_1^T) \rangle_{L^2\times H^{-1}} 
		\\
		- 
		\int_\Omega A \nabla u_0^T \cdot \nabla \psi(x, D_x)^\star  y_0^T \dx
		- 
		\int_\Omega u_1^T \psi(x, D_x)^\star  y_1^T \dx, 
	\end{multline*}
	It is then convenient to notice that the solution $y$ of \eqref{waveequation-control} corresponding to a control function $v$ satisfies, for all $(u_0^T, u_1^T) \in H^1_0(\Omega) \times L^2(\Omega)$,
	$$
		0 
		=
		\int_0^T \int_\omega v \partial_t u \dx\dt
		- 
		\int_\Omega A \nabla u_0^T \cdot \nabla y(T, \cdot) \dx
		-
		\int_\Omega u_1 \partial_t y(T, \cdot) \dx. 
	$$
	Therefore, setting $v = \partial_t U{\large|}_{(0,T) \times \omega}$, the corresponding solution $y$ of \eqref{waveequation-control} satisfies:
	\begin{align}
		& - \div( A \nabla (y(T, \cdot) -\psi(x,D_x)^\star y_0^T)) = - (I - \psi(x, D_x))^\star (I - \psi(x, D_x)) U_0^T, & \text{in }\Omega,
		\label{Formula-yT}
		\\
		&	
		\partial_t y(T, \cdot) - \psi(x, D_x)^\star y_1 = - (I - \psi(x,D_x))^\star (-\Delta_A)^{-1}(I - \psi(x,D_x)) U_1^T, & \text{in } \Omega.
		\label{Formula-dt-yT}
	\end{align}
	
	Accordingly, by elliptic regularity, $ (y(T, \cdot) -\psi(x,D)^\star y_0^T) \in H^2 \cap H^1_0(\Omega)$ and $\partial_t y(0, \cdot) - \psi(x, D)^\star y_1 \in H^1_0(\Omega)$ and we get:
	\begin{align*}
		& \| y(T, \cdot) -\psi(x,D)^\star y_0^T \|_{H^2 \cap H^1_0(\Omega)}
		\leq 
		C \| U_0^T \|_{L^2} \leq C \| (U_0^T, U_1^T)\|_{obs} \leq C \|(y_0^T, y_1^T)\|_{H^1_0\times L^2},
		\\
		& \| \partial_t y(T, \cdot) -\psi(x,D)^\star y_1^T \|_{H^1_0(\Omega)}
		\leq 
		C \| U_1^T \|_{H^{-1}} \leq C \| (U_0^T, U_1^T)\|_{obs} \leq C \|(y_0^T, y_1^T)\|_{H^1_0\times L^2}.
	\end{align*}
	This concludes the proof of Theorem \ref{Thm-Control-With-UC}.
\end{proof} 

\begin{rema}
	In the above proof of Theorem \ref{Thm-Control-With-UC}, we use the duality between the observability and controllability with respect to the pivot space $H^1_0(\Omega) \times L^2(\Omega)$ instead of the usual one developed in \cite{Lions,LionsSIAM88} with respect to the pivot space $L^2(\Omega)$ that we were using in the proof of Theorem \ref{Thm-Obs-local}. This is indeed slightly simpler to handle in the proof of Theorem \ref{Thm-Control-With-UC} since it involves less singular spaces. 
	
One may wonder why this approach was not used in the proof of Theorem \ref{Thm-Obs-local}. The reason lies in the structure of formulas \eqref{Formula-yT}–\eqref{Formula-dt-yT}, which involve commutators with the operator $-\Delta_A$. While these commutators do not affect the regularity of the solutions, they significantly alter their support properties—particularly in the case of formula \eqref{Formula-yT}. As a result, this method is not well-suited for establishing Theorem \ref{Thm-Obs-local}.
\end{rema}

\section{Extensions, open problems and perspectives}

\subsection{Time-dependent coefficients} It would be interesting to investigate the extension of our results to wave equations with time-dependent coefficients. Under a suitable reformulation of the microlocal geometric condition on the pair $(\omega, \mathcal{O})$, the high-frequency propagation results remain valid, and relaxed observability inequalities, similar to those in Lemma \ref{relaxed-obs}, can still be established.

However, removing the compact remainder term in this setting requires a unique continuation result. Assuming analyticity with respect to the time variable, one can obtain a refined observability inequality under a time condition analogous to \eqref{Uniqueness-Condition}. Nevertheless, the compactness-uniqueness argument used in Lemma \ref{UC} is no longer applicable, as the wave equation with time-dependent coefficients is not invariant under time differentiation.

As a result, obtaining sharp observability results analogous to Theorem \ref{Theo-Obs-1} becomes significantly more challenging in the time-dependent case. This limitation is particularly critical when addressing control problems for semilinear or quasilinear wave equations, where time-dependent coefficients naturally arise when applying fixed point techniques.

\subsection{Other observability techniques} Other than the microlocal tools employed in this article, the observability of waves has been often addressed employing multiplier methods \cite{Lions} or Carleman estimates (see for instance \cite{FursikovImanuvilov}). Although they allow to refine global observability estimates when imposing conditions on the support of the initial data (see \cite[Chapitre I, Section 9]{Lions}) by reducing the observability time, these methods do not allow to get the sharp microlocal results in this paper.

\subsection{Schr\"odinger and plate equations} There exists an extensive literature on the observability and control of Schrödinger and plate equations. These models can be roughly viewed as wave-type equations with infinite speed of propagation, which implies that whenever the wave equation is observable or controllable in finite time, the same property holds for the Schrödinger or plate equation in arbitrarily small time, using the same observation and/or control region.

Extending the microlocal and geometric results developed in the present work to such equations remains an interesting and challenging open problem.

\subsection{Control of the heat equation for some specific data} It would be interesting to develop analogues of the results presented in this article for heat-type equations. For instance, one could investigate the cost of controllability in small time for initial data localized in an open subset $\mathcal{O}$, using controls supported on $(0, T) \times \omega$, under the same geometric setting as in Theorem \ref{Theo-Obs-1}. It is natural to conjecture that, in such a case, the controllability cost as $T \to 0$ should be related to the time threshold $T_0$ identified in Theorem \ref{Theo-Obs-1}, and behave like
$C \exp\left( C T_0^2/T \right)$
for some constant $C > 0$.

Indeed, it was shown using the transmutation technique (see \cite{Miller06a}) that one can leverage the controllability properties of the wave equation to derive estimates for the cost of controlling the heat equation in small time. However, the arguments developed in \cite{Miller06a} do not seem directly applicable to the microlocal or geometric setting considered here, and the question remains an open problem.

We also refer to the work \cite{Nguyen-Hoai-Minh-SICON-2022} for a related open question, approached from a different perspective.

Additionally, we note that the transmutation method has also been employed to describe the reachable set for the heat equation (see \cite{Erv-Zua-2011-Transmut}), based on the observability properties of the wave equation. It would be interesting to investigate whether the results of the present work could lead to new estimates on the reachable set for the heat equation, especially in multi-dimensional settings, where this question remains largely open. To our knowledge, the reachable set is fully understood only in the specific case of a ball controlled from its entire boundary, as studied in \cite{Strohmaier-Waters-2022}.

\subsection{Numerical approximation} The numerical analysis of the observability and controllability properties of the wave equation has been also thoroughly investigated. The adaptation of the results in this paper to the discrete context is of interest and would probably require either some suitable filtering processes to avoid the spurious rays (see \cite{Zua05Survey}) and / or some suitable meshes to bend the spurious discrete high-frequency rays (see \cite{Marica-Zuazua-Wigner-2014}).

\subsection{Stabilisation.}It is well known that classical observability and controllability properties are closely linked to the exponential stabilizability of the system. Investigating the stabilization implications of the results developed in this paper thus constitutes an interesting and promising direction for future research.
\bigskip\\

\noindent {\it Acknowledgements.}
B. D was partially supported by the Tunisian Ministry for Higher Education and Scientific Research  within the LR-99-ES20 program.

S. E. is partially supported by the ANR projects TRECOS ANR 20-CE40-0009,  NumOpTes ANR-22-CE46-0005, CHAT ANR-24-CE40-5470.

E. Z.  was funded by then  ERC Advanced Grant CoDeFeL,  the Grants PID2020-112617GB-C22 KiLearn and TED2021-131390B-I00-DasEl of MINECO and PID2023-146872OB-I00-DyCMaMod of MICIU (Spain),  the Alexander von Humboldt-Professorship program, the European Union's Horizon Europe MSCA project ModConFlex, the Transregio 154 Project ``Mathematical Modelling, Simulation and Optimization Using the Example of Gas Networks'' of the DFG, the AFOSR 24IOE027 project, and the Madrid Government - UAM Agreement for the Excellence of the University Research Staff in the context of the V PRICIT (Regional Programme of Research and Technological Innovation).

\end{document}